\documentclass[11pt,a4paper]{article}

\usepackage{url,amsmath,amssymb,latexsym,pstricks,mathrsfs,comment,amsthm,graphicx,tikz,tikz-cd,enumerate,accents,pgffor,cite,wrapfig,multicol,float,cases,calc,bibspacing,geometry,bbm,arydshln}
\usepackage[colorlinks]{hyperref}
\usepackage[T1]{fontenc}
\usepackage{ wasysym }
\usepackage{stmaryrd}
\usepackage{algcompatible}
\usepackage{algorithmicx}
\usepackage{algorithm}
\usepackage[noend]{algpseudocode}
\algrenewcommand{\algorithmiccomment}[1]{\hfill[{\it #1}]}

\usetikzlibrary{calc}

\allowdisplaybreaks

\makeatletter
\newenvironment{myalign}{%
  \def\align@preamble{
    &\hfil
     \strut@
     \setboxz@h{\@lign$\m@th\displaystyle{####}$}%
     \ifmeasuring@\savefieldlength@\fi
     \set@field
     \hfil
     \tabskip\z@skip
    &\setboxz@h{\@lign$\m@th\displaystyle{{}####}$}%
     \ifmeasuring@\savefieldlength@\fi
     \set@field
     \hfil
     \tabskip\alignsep@
  }\align}{\endalign}

\newcommand{\nc}{\newcommand}
\nc{\rnc}{\renewcommand}

\let\oldproofname=\proofname
\rnc{\proofname}{\rm\bf{\oldproofname}}

\rnc{\arraystretch}{1.2}

\rnc{\L}{\mathrel{\mathscr{L}}}
\rnc{\H}{\mathrel{\mathscr{H}}}
\nc{\D}{\mathrel{\mathscr{D}}}
\nc{\R}{\mathrel{\mathscr{R}}}
\nc{\gJ}{\mathrel{\mathscr{J}}}

\nc{\Proj}{\operatorname{Proj}}

\nc{\J}{\mathcal J}
\nc{\M}{\mathcal M}
\nc{\K}{\mathcal K}

\nc{\bn}{\mathbf{n}}

\nc{\bd}{\mathbf{d}}
\nc{\bm}{\mathbf{m}}
\nc{\bo}{\mathbf{O}}
\nc{\bp}{\mathbf{p}}

\nc{\col}[2]{\textstyle{\left[{#1\atop #2}\right]}}

\newcommand{\PP}{\mathscr{P}\mathcal{P}}
\nc{\C}{\mathscr C}
\rnc{\S}{\mathcal S} 
\nc{\B}{\mathcal B} 
\rnc{\P}{\mathcal P} 
\nc{\PB}{\mathcal{PB}} 
\nc{\N}{\mathbb N}
\nc{\Z}{\mathbb Z}
\nc{\F}{\mathbb F}
\rnc{\th}{\theta}
\nc{\al}{\alpha}
\nc{\be}{\beta}
\nc{\Ga}{\Gamma}
\nc{\De}{\Delta}
\nc{\Th}{\Theta}
\nc{\dom}{\operatorname{dom}} 
\nc{\codom}{\operatorname{codom}}
\nc{\coker}{\operatorname{coker}}
\nc{\rank}{\operatorname{rank}}
\nc{\set}[2]{\{ {#1} : {#2} \}} 
\nc{\bigset}[2]{\big\{ {#1}: {#2} \big\}} 
\rnc{\emptyset}{\varnothing}
\nc{\sub}{\subseteq}
\nc{\OR}{\qquad\text{or}\qquad}
\nc{\COMMA}{,\quad}
\nc{\COMMa}{,\ }
\nc{\AND}{\qquad\text{and}\qquad}
\nc{\ANd}{\quad\text{and}\quad}
\rnc{\iff}{\ \Leftrightarrow\ }
\rnc{\implies}{\ \Rightarrow\ }
\nc{\bit}{\begin{itemize}}
\nc{\eit}{\end{itemize}}
\nc{\itemit}[1]{\item[\emph{(#1)}]}
\nc{\itemnit}[1]{\item[(#1)]}
\nc{\pfitem}[1]{\medskip \noindent (#1).}
\nc{\epfres}{\hfill$\Box$}
\nc{\uv}[1]{\fill (#1,2)circle(.17);}
\nc{\lv}[1]{\fill (#1,0)circle(.17);}
\nc{\uvs}[1]{{\foreach \x in {#1} { \uv{\x}}}}
\nc{\lvs}[1]{{\foreach \x in {#1} { \lv{\x}}}}
\nc{\stline}[2]{\draw(#1,2)--(#2,0);}
\nc{\stlines}[1]{{\foreach \x/\y in {#1} { \stline{\x}{\y} }}}
\nc{\darcx}[3]{\draw(#1,0)arc(180:90:#3) (#1+#3,#3)--(#2-#3,#3) (#2-#3,#3) arc(90:0:#3);}
\nc{\darc}[2]{\darcx{#1}{#2}{.4}}
\nc{\uarcx}[3]{\draw(#1,2)arc(180:270:#3) (#1+#3,2-#3)--(#2-#3,2-#3) (#2-#3,2-#3) arc(270:360:#3);}
\nc{\uarc}[2]{\uarcx{#1}{#2}{.4}}

\newcommand{\darcxred}[3]{\draw[red](#1,0)arc(180:90:#3) (#1+#3,#3)--(#2-#3,#3) (#2-#3,#3) arc(90:0:#3);}
\newcommand{\uarcxred}[3]{\draw[red](#1,2)arc(180:270:#3) (#1+#3,2-#3)--(#2-#3,2-#3) (#2-#3,2-#3) arc(270:360:#3);}
\newcommand{\duarcx}[3]{\draw(#1,0)arc(180:270:#3) (#1+#3,0-#3)--(#2-#3,0-#3) (#2-#3,0-#3) arc(270:360:#3);}
\newcommand{\duarcxred}[3]{\draw[red](#1,0)arc(180:270:#3) (#1+#3,0-#3)--(#2-#3,0-#3) (#2-#3,0-#3) arc(270:360:#3);}
\newcommand{\duarc}[2]{\duarcx{#1}{#2}{.4}}

\newcommand{\stlined}[1]{\draw(#1,0)--(#1,-0.75);}
\newcommand{\stlineu}[1]{\draw(#1,0)--(#1,0.75);}

\newcommand{\bluebox}[4]{
\draw[color=blue!20, fill=blue!20] (#1,#2)--(#3,#2)--(#3,#4)--(#1,#4)--(#1,#2);
}
\newcommand{\bluetrap}[8]{
\draw[color=blue!20, fill=blue!20] (#1,#2)--(#3,#4)--(#5,#6)--(#7,#8)--(#1,#2);
}
\newcommand{\stlinedot}[2]{\draw[dotted](#1,2)--(#2,0);}
\newcommand{\stlinedots}[1]{
{\foreach \x/\y in {#1}
{ \stlinedot{\x}{\y} }
}
}
\newcommand{\uvertlabel}[2]{\draw(#1,2+.5)node{{\tiny $#2$}};}
\newcommand{\uvertlabels}[1]{
{\foreach \x/\y in {#1}
{ \uvertlabel{\x}{\y} }
}
}
\newcommand{\lvertlabel}[2]{\draw(#1,0-.5)node{{\tiny $#2$}};}
\newcommand{\lvertlabels}[1]{
{\foreach \x/\y in {#1}
{ \lvertlabel{\x}{\y} }
}
}
\newcommand{\udotted}[2]{\draw [dotted] (#1,2)--(#2,2);}
\newcommand{\udotteds}[1]{{
\foreach \x/\y in {#1}
{ \udotted{\x}{\y}}
}}
\newcommand{\ldotted}[2]{\draw [dotted] (#1,0)--(#2,0);}
\newcommand{\ldotteds}[1]{{
\foreach \x/\y in {#1}
{ \ldotted{\x}{\y}}
}}

\numberwithin{equation}{section}

\newtheorem{thm}[equation]{Theorem}
\newtheorem{lemma}[equation]{Lemma}
\newtheorem{cor}[equation]{Corollary}
\newtheorem{prop}[equation]{Proposition}

\theoremstyle{definition}

\newtheorem{rem}[equation]{Remark}

\begin{document}

\title{Enumeration of idempotents in planar diagram monoids~\vspace{-5ex}}
\author{}
\date{}

\maketitle
\begin{center}
{\large 
Igor Dolinka,%
\hspace{-.3em}\footnote{Department of Mathematics and Informatics, University of Novi Sad, Trg Dositeja Obradovi\'ca 4, 21101 Novi Sad, Serbia. {\it Email:} {\tt dockie@dmi.uns.ac.rs}}
James East,%
\hspace{-.3em}\footnote{Centre for Research in Mathematics, School of Computing, Engineering and Mathematics, University of Western Sydney, Locked Bag 1797, Penrith NSW 2751, Australia. {\it Email:} {\tt J.East@uws.edu.au}}
Athanasios Evangelou,%
\hspace{-.3em}\footnote{\label{footnote:AE}School of Mathematics and Physics, University of Tasmania, Private Bag 37, Hobart 7001, Australia. {\it Emails:} {\tt aoost@utas.edu.au}, {\tt D.FitzGerald@utas.edu.au}, {\tt nicholas.ham@utas.edu.au}}
Des FitzGerald,\hspace{-.2em}\textsuperscript{\ref{footnote:AE}}}\\
{\large 
Nicholas Ham,\hspace{-.2em}\textsuperscript{\ref{footnote:AE}}
James Hyde,%
\hspace{-.3em}\footnote{\label{footnote:JH}School of Mathematics and Statistics, University of St Andrews, St Andrews, UK. {\it Emails:} {\tt jth4@st-andrews.ac.uk}, {\tt jdm3@st-andrews.ac.uk}}
Nicholas Loughlin,%
\hspace{-.3em}\footnote{School of Mathematics and Statistics, Newcastle University, Newcastle NE1 7RU, UK.  {\it Email:} {\tt n.j.loughlin@newcastle.ac.uk}}
James D.~Mitchell,\hspace{-.2em}\textsuperscript{\ref{footnote:JH}}
}\\
~\\
{\large \today}
\end{center}

\begin{abstract}
We classify and enumerate the idempotents in several planar diagram monoids: namely, the Motzkin, Jones (a.k.a.~Temperley-Lieb) and Kauffman monoids.  The classification is in terms of certain vertex- and edge-coloured graphs associated to Motzkin diagrams.  The enumeration is necessarily algorithmic in nature, and is based on parameters associated to cycle components of these graphs.  We compare our algorithms to existing algorithms for enumerating idempotents in arbitrary (regular $*$-) semigroups, and give several tables of calculated values.


\end{abstract}

\section{Introduction}\label{sect:intro}

Diagram monoids arise in numerous branches of mathematics and science, including representation theory, statistical mechanics and knot theory \cite{HR2005,Brauer1937,Jones1994_2,Martin1994,TL1971, Jones1987,Kauffman1990,LZ2012}.  
Many studies of diagram monoids have been combinatorial in nature \cite{EMRT2018,DEG2017,EG2017,DEEFHHL1,DE5,Maltcev2007,JEpnsn}, and idempotents have played a large role in several of these works.  Historically, it is interesting to note that products of idempotents \cite{Maltcev2007,JEpnsn,EF2012} were understood several years before the idempotents themselves were \cite{DEEFHHL1}; this is largely due to the fact that diagram monoids have natural anti-involutions that give them regular $*$-semigroup structures \cite{EF2012,NS1978}, meaning that arbitrary idempotents are products of simpler idempotents known as projections (see Section \ref{sect-existing} for definitions).

In \cite{DEEFHHL1}, classifications and enumerations were given for the idempotents in the partition, Brauer and partial Brauer monoids, and also for the idempotent basis elements in the corresponding diagram algebras.  The current paper continues in this direction, focusing this time on planar diagram monoids, such as the Motzkin, Jones (a.k.a.~Temperley-Lieb) and Kauffman monoids.
However, the methods employed here are necessarily different to those of \cite{DEEFHHL1}, as the planarity constraint means that the set symmetries used to study the monoids in \cite{DEEFHHL1} are no longer available.  In fact, we believe that enumeration of the planar idempotents is inimical to closed-form solution; for one thing, the highly complex meandric numbers \cite{DFGG1997,LaCroix2003} occur during the enumeration of Kauffman idempotents, as noted below in Section \ref{sect-values}.  We seek to fill this gap by presenting methods for computing the numbers of idempotents in the Motzkin, Jones and Kauffman monoids, with attention given to efficiency of the algorithms involved.  

The article is organised as follows.
In Section \ref{sect-existing}, we describe two existing algorithms for enumerating the idempotents of arbitrary finite (regular $*$-) semigroups, commenting on their respective complexities in general and in the context of the Jones and Motzkin monoids in particular.  
%
Section \ref{sect-theory} develops a theory of idempotents in the Motzkin, Jones and Kauffman monoids.  
A key role is played by certain graphs, called interface graphs, associated to arbitrary Motzkin elements.  These graphs are used to classify the idempotents in Proposition~\ref{prop-idempotents} and Corollaries \ref{cor-jones-idempotents} and \ref{cor-kauffman-idempotents}, and then to enumerate them in Theorems~\ref{thm-EMJ},~\ref{thm-EK} and \ref{thm-ErMJK}; see also Proposition \ref{prop-D-1}.  
Section \ref{sect-algorithms} presents a number of algorithms, based on the theory developed in Section \ref{sect-theory}, that calculate the number of idempotents in the Jones, Kauffman and Motzkin monoids; C++ implementations of these algorithms can be found at \cite{Mitchell2016aa}.  Finally, Section \ref{sect-values} gives several tables of calculated values, including comparative run-times of the various algorithms.

The reader is referred to the monographs of Higgins \cite{Higgins} and Howie \cite{Howie} for background on semigroups in general, and to the introduction of \cite{DEEFHHL1} and references therein---in particular to foundational articles of Jones \cite{Jones1994_2}, Martin \cite{Martin1994} and Halverson and Ram \cite{HR2005}---for background and relevant detail on the partition, Brauer and partial Brauer monoids.

\section{Existing algorithms}\label{sect-existing}

In this section, we describe two existing approaches to counting idempotents in semigroups such as those we study in this article.  The first applies to any finite semigroup, while the second applies to any finite regular $*$-semigroup; see below for the definitions.  
In Section \ref{sect-values}, we will discuss the performance of these two approaches when applied to the diagram monoids we are concerned with, and we will compare them with the new algorithms presented in Section \ref{sect-algorithms}. 

The first approach to counting idempotents in an arbitrary finite semigroup $S$ is simple: create the elements of $S$, and then check if $x ^ 2 = x$ for each $x \in S$; see Algorithm~\ref{algorithm-basic}. If the semigroup $S$ is generated by
$A\subseteq S$, then $S$ can be enumerated using the Froidure-Pin
Algorithm~\cite{Froidure1997aa}. If the complexity of
multiplying elements in $S$ is assumed to be constant, then the complexity of the Froidure-Pin
Algorithm is $O\big(|S||A|\big)$, so
Algorithm~\ref{algorithm-basic} has complexity $O\big(|S||A|+|S|\big)=O\big(|S||A|\big)$.

\begin{algorithm}[h]
  \caption{Count idempotents in a semigroup $S$}
  \label{algorithm-basic}
  \begin{algorithmic}[1]
    \State $n := 0$
    \For{$x\in S$}
    \If{$x ^ 2 = x$}
    \State $n\gets n + 1$
    \EndIf
    \EndFor
    \State \Return $n$
  \end{algorithmic}
\end{algorithm}

The approach just described requires the creation of each element of $S$ in order to check which elements are idempotents.  When $S$ is very large, this can be impractical, in terms of both space and time.  The second approach improves on the first in the case that $S$ is a regular $*$-semigroup.  To describe it, we must first recall some background.

Recall from \cite{NS1978} that a semigroup $S$ is a \emph{regular $*$-semigroup} if there is a unary operation ${}^*:S\to S$ such that $x^{**}=x$, $(xy)^*=y^*x^*$ and $xx^*x=x$ for all $x,y\in S$.  For the remainder of this section, we fix a finite regular $*$-semigroup $S$.  Recall that \emph{Green's relations} $\R$, $\L$, $\gJ$, $\H$ and $\D$ are defined on~$S$ as follows.  Let $x,y\in S$.  We say that $x\R y$ if $xS=yS$, that $x\L y$ if $Sx=Sy$, and that $x\gJ y$ if $SxS=SyS$.  The relation $\H$ is defined to be the intersection of $\R$ and $\L$, while $\D$ is defined to be the join of $\R$ and $\L$; that is, $\D$ is the least equivalence on $S$ containing both $\R$ and $\L$.  It is well known that ${\D}={\R}\circ{\L}={\L}\circ{\R}$, and that ${\D}={\gJ}$ since $S$ is finite; see \cite[Chapter 2]{Howie} for more background on Green's relations.  
We write $E(S)=\set{x\in S}{x^2=x}$ for the set of all idempotents of $S$.
An element $x\in S$ is called a \emph{projection} if $x^2=x=x^*$.  The set of all projections of $S$ is denoted by $\Proj(S)$.  
The next result is true of any regular $*$-semigroup, whether finite or infinite; proofs of the various parts may be found in \cite{Imaoka1980,NS1978,Petrich1985}.

\begin{lemma}\label{lem-RSS}
Let $S$ be a regular $*$-semigroup.  Then
\bit
\itemit{i} $\Proj(S)=\set{xx^*}{x\in S}=\set{x^*x}{x\in S}$,
\itemit{ii} $E(S)=\set{xy}{x,y\in\Proj(S)}$,
\itemit{iii} every element of $S$ is $\R$-related to a unique projection,
\itemit{iv} every element of $S$ is $\L$-related to a unique projection,
\itemit{v} for any $x\in\Proj(S)$ and any $a\in S$, $a^*xa\in\Proj(S)$,
\itemit{vi} for any $x,y\in S$, $x\R y$ if and only if $xx^*=yy^*$,
\itemit{vii} for any $x,y\in S$, $x\L y$ if and only if $x^*x=y^*y$. \epfres
\eit
\end{lemma}

It follows from Lemma \ref{lem-RSS}(v), and the identity $(ab)^*=b^*a^*$, that $S$ has a right action on $\Proj(S)$, defined by
\begin{equation}\label{equation-action}
  x\cdot a = a^*xa  \qquad\text{for $x\in\Proj(S)$ and $a\in S$.}
\end{equation}
For $x\in\Proj(S)$, we will write 
\[
[x] = \set{y\in\Proj(S)}{\text{$x=y\cdot a$ and $y=x\cdot b$ for some $a,b\in S$}}
\]
for the \emph{strongly connected component} of $x$ under the action~\eqref{equation-action}.
For any subset $A\sub S$, we will write $E(A)=A\cap E(S)$ and $\Proj(A)=A\cap\Proj(S)$ for the set of all idempotents and projections belonging to~$A$, respectively.  If $\mathscr K$ is any of Green's relations, and if $x\in S$, we denote the $\mathscr K$-class of $x$ by $K_x$. 
Since $S$ is finite, it follows that $S$ has the \emph{stability} property: namely, for any $x,y\in S$, $xy\D x$ implies $xy\R x$, and $xy\D y$ implies $xy\L y$; see \cite[Section A.2]{RSbook}.  The various parts of the next lemma may be known, but we give proofs for completeness.

\begin{lemma}\label{lem-RSS2}
Let $S$ be a finite regular $*$-semigroup.  Then
\bit
\itemit{i} for any $x\in\Proj(S)$, $[x] = \Proj(D_x) = \set{y\in\Proj(S)}{x\D y}$,
\itemit{ii} for any $\D$-class $D$ of $S$, and for any projections $x,y\in\Proj(D)$, $R_x\cap L_y$ contains an idempotent if and only if $xy\D x$, in which case this idempotent is $xy$,
\itemit{iii} for any $\D$-class $D$ of $S$, the number of idempotents in $D$ is equal to the cardinality of the set
\[
{\set{(x,y)}{x,y\in\Proj(D),\ xy\D x}}.
\]
\eit
\end{lemma}

\begin{proof}
(i).
Let $x\in\Proj(S)$.  If $y\in[x]$, then $x=y\cdot a=a^*ya$ and $y=x\cdot b=b^*xb$ for some $a,b\in S$, so that $x\gJ y$, whence $x\D y$ as $S$ is finite.  Conversely, suppose $y\in\Proj(S)$ and $x\D y$.  Then $x\R a$ and $a\L y$ for some $a\in S$.  By Lemma \ref{lem-RSS}(vi) and (vii), and since $x,y\in\Proj(S)$, we obtain $aa^*=xx^*=x$ and $a^*a=y^*y=y$.  It follows that $y=a^*a=a^*aa^*a=a^*xa=x\cdot a$, and similarly $x=y\cdot a^*$, so that $y\in[x]$.

\pfitem{ii}
Suppose first that $R_x\cap L_y$ contains some idempotent $e$.  Since $e\R x$, Lemma \ref{lem-RSS}(vi) gives $ee^*=xx^*=x$, and similarly $e^*e=y$.  Then $e=ee^*e=e(ee)^*e=(ee^*)(e^*e)=xy$.  Also, $xy=e\R x$ implies $xy\D x$.

Conversely, suppose $xy\D x$.  Then stability gives $xy\R x$.  Since $xy\D x\D y$, stability also gives $xy\L y$.  Thus, $xy\in R_x\cap L_y$.  Lemma \ref{lem-RSS}(ii) gives $xy\in E(S)$.

\pfitem{iii}
First note that the number of idempotents in $D$ is equal to the number of $\H$-classes in~$D$ containing an idempotent, since each $\H$-class contains at most one idempotent.  If $a\in D$, then $H_a=R_a\cap L_a$, and by Lemma \ref{lem-RSS}(iii) and~(iv) we have $R_a=R_x$ and $L_a=L_y$ for unique projections $x,y\in\Proj(D)$.  That is, every $\H$-class in $D$ is equal to $R_x\cap L_y$ for unique projections $x,y\in\Proj(D)$.  By part (ii), just proved, this $\H$-class contains an idempotent if and only if $xy\D x$. 
\end{proof}

Parts (i) and (iii) of Lemma \ref{lem-RSS2} form the basis of the second approach to computing the number of idempotents in the finite regular $*$-semigroup $S$, given by a generating set $A\sub S$; see Algorithm \ref{algorithm-regular-*}.  
Roughly speaking, the steps of this algorithm are:
\bit
\itemnit{1}  We first create $\Proj(S)$.  
\itemnit{2}  We then create the sets $\Proj(D)$, as $D$ runs over the set of all $\D$-classes of~$S$.  
\itemnit{3}  For each $\D$-class $D$, we then find the cardinality of the set given in Lemma~\ref{lem-RSS2}(iii), and sum over all~$D$.  
\eit

Step (1) can be achieved using the action from~\eqref{equation-action} in a simple orbit algorithm whose input is the
generators~$A$; see \cite[Algorithm 1]{East2015aa}.  If the complexity of determining $x\cdot a$ is assumed to be constant, then the complexity of \cite[Algorithm 1]{East2015aa}, and hence the complexity of Step (1), is $O\big(|{\Proj(S)}||A|\big)$.

By Lemma~\ref{lem-RSS2}(i), the sets $\Proj(D)$ correspond to the strongly connected components of the action of~$S$ on $\Proj(S)$ given in~\eqref{equation-action}.  These can be found using standard algorithms from graph theory, such as Tarjan's \cite{Tarjan1972aa} or Gabow's \cite{Gabow2000aa}, for example; see also the monograph of Sedgewick \cite{Sedgewick1983aa}. The complexity of these algorithms, and thus the complexity of Step (2), is $O\big(|{\Proj(S)}| + |A|\big)$, which is bounded above by $O\big(|{\Proj(S)}||A|\big)$, the complexity of Step (1).

If the $\D$-classes of $S$ are $D_1,\ldots,D_r$, and if these $\D$-classes have $m_1,\ldots,m_r$ projections, respectively, then Step (3) involves $m_1^2+\cdots+m_r^2$ products and checks for $\D$-relatedness (modulo some optimisations discussed below).  Thus, the total complexity of Algorithm \ref{algorithm-regular-*} is
\[
O\big(|{\Proj(S)}||A| + m_1^2+\cdots+m_r^2\big).
\]
The $\H$-classes in a single $\D$-class of a semigroup all have the same size \cite[Lemma 2.2.3]{Howie}.  If the $\H$-classes in the $\D$-class $D_i$ have size $h_i$, then $|D_i|=m_i^2h_i$, since $D_i$ has $m_i$ $\R$- and $\L$-classes (Lemma \ref{lem-RSS}(iii) and~(iv)), and so
\[
m_1^2+\cdots+m_r^2\leq m_1^2h_1+\cdots+m_r^2h_r=|D_1|+\cdots+|D_r|=|S|.
\]
This upper bound is realised if and only if $h_i=1$ for all $i$: i.e., if $S$ is $\H$-trivial.  In this worst case, the total complexity of Algorithm \ref{algorithm-regular-*} is $O\big(|{\Proj(S)}||A| + |S|\big)$.  When we compare this to the complexity of Algorithm~\ref{algorithm-basic}, which we noted above was $O\big(|S||A|\big)$, we see that Algorithm~\ref{algorithm-regular-*} has a significant advantage if $|{\Proj(S)}|$ is small relative to $|S|$.  Note that $|{\Proj(S)}|=m_1+\cdots+m_r$.

We note that Algorithm \ref{algorithm-regular-*}, presented below, contains a number of simple optimisations.  First, the $\H$-class $H_x=R_x\cap L_x$ of a projection $x\in\Proj(D)$ always contains an idempotent: namely, $x$ itself (see Line~\ref{algorithm-regular-*-line6}).  Secondly, if $x,y\in\Proj(D)$, then $xy\in D$ if and only if $yx=y^*x^*=(xy)^*\in D$, so we only need to test one of $xy$ or $yx$ for membership in $D$ (see Lines \ref{algorithm-regular-*-line7}--\ref{algorithm-regular-*-line10}).  We also note that when $S$ is any of the diagram monoids we consider, the $\D$-relation is given by equality of the \emph{ranks} of elements of $S$, and is easily checked computationally; see Section \ref{sect-theory} for the definition of rank, and also \cite{Wilcox2007,LF2006,DEG2017} for more on Green's relations on diagram monoids.
Finally, we note that Algorithm~\ref{algorithm-regular-*} can be derived from \cite[Algorithm 10]{East2015aa}, which counts idempotents in a fixed $\R$-class.

\begin{algorithm}[H]
  \caption{Count idempotents in a regular $*$-semigroup $S$}
  \label{algorithm-regular-*}
  \begin{algorithmic}[1]
    \State Find $\Proj(S)$ 
    \State Find the strongly connected components $C_1, \ldots, C_r$ of the action of $S$ on     $\Proj(S)$ defined in \eqref{equation-action}
    \State $n := 0$
    \For{$i \in \{1, \ldots, r\}$}
    \State if $C_i = \{x_1, x_2, \ldots, x_m\}$
    \State $n \gets n+m$   \label{algorithm-regular-*-line6}
    \For{$j \in\{1, \ldots, m\}$}   \label{algorithm-regular-*-line7}
    \For{$k \in \{j + 1, \ldots, m\}$}
    \If{$x_ix_j \D x_i$}
    \State $n \gets n + 2$   \label{algorithm-regular-*-line10}
    \EndIf
    \EndFor
    \EndFor
    \EndFor
    \State \Return $n$
  \end{algorithmic}
\end{algorithm}

This paper mostly concerns the case in which $S$ is a Jones, Motzkin or Kauffman monoid.  As noted above, the definitions of these monoids are given in Section \ref{sect-theory}, but here we make some brief comments relevant to the current discussion.  For each non-negative integer $n$, we have a regular $*$-monoid $\J_n$ (Jones) and $\M_n$ (Motzkin).
The sizes of these monoids, and the sizes of their sets of projections, are given (see \cite{DEG2017,BH2014,EG2017}) by
\[
|\J_n| = C_n \COMMA |\M_n|=\mu(2n,0) \COMMA |{\Proj}(\J_n)| = \sum_{r=0}^n \frac{r+1}{n+1} \binom{n+1}{\frac{n-r}2} \COMMA |{\Proj}(\M_n)| = \sum_{r=0}^n \mu(n,r).
\]
Here, $C_n=\frac1{n+1}\binom{2n}n$ is the $n$th Catalan number; we interpret a binomial coefficient $\binom mk$ to be $0$ if $k$ is not an integer between $0$ and $m$; and the \emph{Motzkin triangle} numbers $\mu(n,r)$ are defined by the recurrence
\begin{myalign}
\nonumber \mu(0,0)=1 \COMMA
\mu(n,r)=0 &\qquad\qquad\text{if $n<r$ or $r<0$,}\\
\label{equation-motzkin1} \mu(n,r)=\mu(n-1,r-1)+\mu(n-1,r)+\mu(n-1,r+1) &\qquad\qquad\text{if $0\leq r\leq n$ and $n\geq1$.}
\end{myalign}
The numbers $\mu(n,r)$ are also given by the formula
\begin{equation}\label{equation-motzkin2}
\mu(n, r) = \sum_{j = 0} ^ n \binom nj \left[ \binom{n-j}{r+j}-\binom{n-j}{r+j+2} \right].
\end{equation}
See Sequences A000108, A001006, A026300 in \cite{OEIS}.  
%
Values of $|S|$ and $|{\Proj}(S)|$ are given in Table~\ref{tab-comp-1} for $S=\J_n$ or $\M_n$ with $n\leq15$, by way of indicating the relative complexities of Algorithms~\ref{algorithm-basic} and~\ref{algorithm-regular-*}.
Run times are given in Section \ref{sect-values}, further highlighting the advantage of Algorithm \ref{algorithm-regular-*} over Algorithm \ref{algorithm-basic} in these cases.

\begin{table}[h]
  \begin{center}
\begin{tabular}{|r|r|r|r|r|}
\hline
  \multicolumn{1}{|c|}{$n$}   & \multicolumn{1}{c|}{$|\J_n|$}     & \multicolumn{1}{c|}{$|{\Proj}(\J_n)|$} & \multicolumn{1}{c|}{$|\M_n|$} & \multicolumn{1}{c|}{$|{\Proj}(\M_n)|$} \\ \hline
0    &  1          &  1             &  1              &1          \\ 
1    &  1          &  1             &  2              &1          \\ 
2    &  2          &  2             &  9              &2          \\  
3    &  5          &  3             &  51             &5          \\ 
   4    &  14         &  6          &  323            &13         \\
   5    &  42         &  10         &  2188           &35         \\
\hline
   6    &  132        &  20         &  15\ 511          &96          \\
   7    &  429        &  35         &  113\ 634         &267          \\
   8    &  1430       &  70         &  853\ 467         &750          \\
   9    &  4862       &  126        &  6536\ 382        &2123          \\
   10   &  16\ 796      &  252        &  50\ 852\ 019       &6046           \\
\hline
   11   &  58\ 786      &  462        &  400\ 763\ 223      &17\ 303           \\
   12   &  208\ 012     &  924        &  3192\ 727\ 797     &49\ 721            \\
   13   &  742\ 900     &  1716       &  25\ 669\ 818\ 476    &143\ 365            \\
   14   &  2674\ 440    &  3432       &  208\ 023\ 278\ 209   &414\ 584             \\
   15   &  9694\ 845    &  6435       &  1697\ 385\ 471\ 211  &1201\ 917             \\
\hline
\end{tabular}                          
\end{center}
\vspace{-5mm}
  \caption{The sizes of the Jones and Motzkin monoids, $\J_n$ and $\M_n$, and of their sets of projections, $\Proj(\J_n)$ and $\Proj(\M_n)$.}
  \label{tab-comp-1}
\end{table}

Algorithms~\ref{algorithm-basic} and~\ref{algorithm-regular-*}, as presented above, are both
\emph{embarrassingly parallel}.  Parallel versions of the algorithms are
implemented in the Semigroups package for GAP \cite{GAP}.

Finally, we note that the Kauffman monoid $\K_n$ (also defined in Section \ref{sect-theory}) is infinite, and also not a regular $*$-semigroup, so neither of the algorithms discussed in this section apply to it.  It is possible to define a finite quotient of~$\K_n$ that has only one more idempotent than $\K_n$ \cite{LF2006}, but this quotient is still not a regular $*$-semigroup, so only Algorithm~\ref{algorithm-basic} would apply.

\section{Idempotents of planar diagram monoids}\label{sect-theory}

In this section, we define the diagram monoids we will be concerned with, before describing methods to classify and enumerate the idempotents of these monoids.  The classification involves certain graphs, called \emph{interface graphs}, associated to arbitrary Motzkin elements.  The enumeration is based on natural parameters associated to certain cycle components of the interface graphs, as well as a map that sends Motzkin idempotents to lower-rank idempotents.

\subsection{Definitions and preliminaries}

Let $n$ be a positive integer, and write $\bn=\{1,\ldots,n\}$ and $\bn'=\{1',\ldots,n'\}$.  The \emph{partition monoid of degree~$n$}, denoted $\P_n$, is the monoid of all set partitions of $\bn\cup\bn'$ under a product described below.  Thus, an element of $\P_n$ is a set $\alpha=\{A_1,\ldots,A_k\}$, for some $k$, where the $A_i$ are pairwise disjoint non-empty subsets of $\bn\cup\bn'$ whose union is all of $\bn\cup\bn'$; the $A_i$ are called the \emph{blocks} of $\alpha$.  By convention, $\P_0$ contains a single element, the empty partition, but we will assume $n\geq1$ since all results concerning $\P_0$ are trivial.

\begin{figure}[h]
\begin{center}
\begin{tikzpicture}[scale=.5]

\begin{scope}[shift={(0,0)}]	
\uvs{1,...,6}
\lvs{1,...,6}
\uarcx14{.6}
\uarcx23{.3}
\uarcx56{.3}
\darc12
\darcx26{.6}
\darcx45{.3}
\stline34
\draw(0.6,1)node[left]{$\alpha=$};
\draw[->](7.5,-1)--(9.5,-1);
\end{scope}

\begin{scope}[shift={(0,-4)}]	
\uvs{1,...,6}
\lvs{1,...,6}
\uarc12
\uarc34
\darc45
\darc56
\stline31
\stline55
\draw(0.6,1)node[left]{$\beta=$};
\end{scope}

\begin{scope}[shift={(10,-1)}]	
\uvs{1,...,6}
\lvs{1,...,6}
\uarcx14{.6}
\uarcx23{.3}
\uarcx56{.3}
\darc12
\darcx26{.6}
\darcx45{.3}
\stline34
\draw[->](7.5,0)--(9.5,0);
\end{scope}

\begin{scope}[shift={(10,-3)}]	
\uvs{1,...,6}
\lvs{1,...,6}
\uarc12
\uarc34
\darc45
\darc56
\stline31
\stline55
\end{scope}

\begin{scope}[shift={(20,-2)}]	
\uvs{1,...,6}
\lvs{1,...,6}
\uarcx14{.6}
\uarcx23{.3}
\uarcx56{.3}
\darc14
\darc45
\darc56
\stline21
\draw(6.4,1)node[right]{$=\alpha\beta$};
\end{scope}

\end{tikzpicture}
\end{center}
\vspace{-5mm}
\caption{Two partitions $\alpha,\beta\in\P_6$ (left), their product $\alpha\beta\in\P_6$ (right), and the product graph $\Pi(\alpha,\beta)$ (centre).}
\label{fig:P6}
\end{figure}
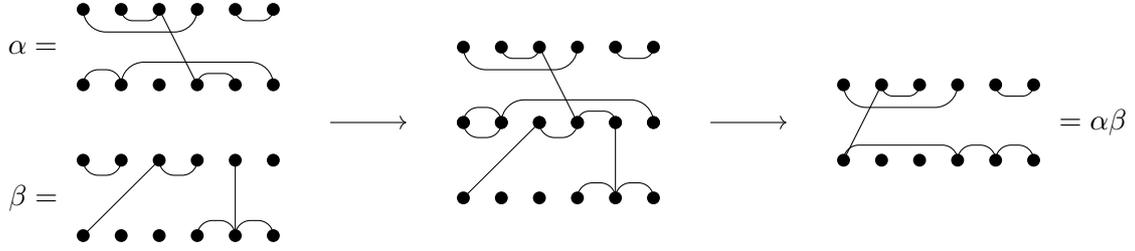

A partition $\alpha\in\P_n$ may be pictured (non-uniquely) as a graph with vertex set $\bn\cup\bn'$, and with any edge set having the property that the connected components of the graph correspond to the blocks of the partition; the vertices of such a graph are always drawn with $1,\ldots,n$ on an upper row, increasing from left to right, and vertices $1',\ldots,n'$ directly below.   For example, the partitions
\begin{align*}
\alpha = \big\{ \{1,4\},\{2,3,4',5'\},\{5,6\},\{1',2',6'\},\{3'\}\big\} \ANd
\beta = \big\{ \{1,2\}, \{3,4,1'\}, \{5,4',5',6'\}, \{6\}, \{2'\}, \{3'\} \big\}
\end{align*}
from $\P_6$ are pictured in Figure \ref{fig:P6}.  As usual, we will generally identify a partition with any graph representing it.

A block $A$ of a partition $\alpha$ is referred to as a \emph{transversal} if $A\cap\bn\not=\emptyset$ and $A\cap\bn'\not=\emptyset$, or a \emph{non-transversal} otherwise.  
For example, $\alpha\in\P_6$ defined above has $\{2,3,4',5'\}$ as its only transversal, and has upper non-transversals $\{1,4\}$ and $\{5,6\}$, and lower non-transversals $\{1',2',6'\}$ and~$\{3'\}$.

The \emph{domain} and \emph{codomain} of $\alpha\in\P_n$ are the subsets of $\bn$ defined by
\begin{align*}
\dom(\alpha) &= \set{i\in\bn}{\text{$i$ belongs to a transversal of $\alpha$}}, \\
\codom(\alpha) &= \set{i\in\bn}{\text{$i'$ belongs to a transversal of $\alpha$}} .
\end{align*}
The \emph{rank} of $\alpha\in\P_n$, denoted $\rank(\alpha)$, is defined to be the number of transversals of $\alpha$.  For example, with $\al\in\P_6$ as defined above,  $\rank(\al)=1$, $\dom(\alpha)=\{2,3\}$ and $\codom(\alpha) = \{4,5\}$.

The product of two partitions $\alpha,\beta\in\P_n$ is defined as follows.  
Write $\bn''=\{1'',\ldots,n''\}$.  Let $\alpha^\vee$ be the graph obtained from $\alpha$ by changing the label of each lower vertex $i'$ to~$i''$, and let $\beta^\wedge$ be the graph obtained from $\beta$ by changing the label of each upper vertex~$i$ to~$i''$.  Consider now the graph $\Pi(\alpha,\beta)$ on the vertex set~$\bn\cup \bn'\cup \bn''$ obtained by joining $\alpha^\vee$ and~$\beta^\wedge$ together so that each lower vertex $i''$ of $\alpha^\vee$ is identified with the corresponding upper vertex $i''$ of $\beta^\wedge$.  Note that $\Pi(\alpha,\beta)$, which we call the \emph{product graph}, may contain pairs of parallel edges.  We define $\alpha\beta\in\P_n$ to be the partition satisfying the property that $x,y\in\bn\cup\bn'$ belong to the same block of $\alpha\beta$ if and only if $x$ and $y$ are connected by a path in $\Pi(\alpha,\beta)$.  
This process is illustrated in Figure~\ref{fig:P6}, with $\alpha,\beta\in\P_6$ defined above.
The operation is associative, so~$\P_n$ is a semigroup: in fact, a monoid, with identity element $\big\{\{1,1'\},\ldots,\{n,n'\}\big\}$.  

Note that the product graph $\Pi(\alpha,\beta)$ may contain connected components that only involve vertices from~$\bn''$; these are called \emph{floating components} of $\Pi(\alpha,\beta)$.  These play a crucial role in the definition of the \emph{partition algebras}, and also the \emph{twisted partition monoids}, which we now describe.  Specifically, if $\tau(\alpha,\beta)$ denotes the number of floating components in the product graph $\Pi(\alpha,\beta)$, then one easily checks that
\[
\tau(\alpha,\beta)+\tau(\alpha\beta,\gamma)=\tau(\alpha,\beta\gamma)+\tau(\beta,\gamma) \qquad\text{for all $\alpha,\beta,\gamma\in\P_n$.}
\]
It then follows that the product $\star$ defined on the set $\P_n^\tau=\N\times\P_n=\set{(i,\alpha)}{i\in\N,\ \alpha\in\P_n}$ by
\begin{equation}\label{eq:star}
(i,\alpha)\star(j,\beta)=\big(i+j+\tau(\alpha,\beta),\alpha\beta\big)
\end{equation}
is associative.  (Here, $\N=\{0,1,2,\ldots\}$ denotes the set of natural numbers.)  We call $\P_n^\tau$ with this product the \emph{twisted partition monoid}.  

A partition $\alpha\in\P_n$ is \emph{planar} if there is a graphical representation of $\alpha$ in which:
\begin{itemize}
\item[(i)] all the edges are drawn within the rectangle determined by the vertices; and
\item[(ii)] there are no crossings within the interior of the rectangle.
\end{itemize}
For example, of the two partitions $\alpha,\beta\in\P_6$ defined above, $\alpha$ is not planar, but $\beta$ is.  Since the product of two planar partitions is clearly planar, it follows that the set of all such planar partitions forms a submonoid of $\P_n$, and we denote this planar submonoid by $\PP_n$.

The \emph{partial Brauer monoid} and the \emph{Brauer monoid} are the submonoids of $\P_n$ defined by
\begin{align*}
\PB_n = \set{\alpha\in\P_n}{\text{all blocks of $\alpha$ have size $1$ or $2$}} &\AND
\B_n = \set{\alpha\in\P_n}{\text{all blocks of $\alpha$ have size $2$}}.
\intertext{The \emph{Motzkin monoid} and the \emph{Jones monoid} are the planar submonoids of $\P_n$ defined by}
\M_n = \PB_n\cap\PP_n &\AND \J_n = \B_n\cap\PP_n.
\end{align*}
A Motzkin element $\gamma\in\M_{20}$ is pictured in Figure \ref{fig:M20}; the reader is invited to verify that $\gamma$ is in fact an idempotent. 
It is well known that~$\PP_n$ is isomorphic to $\J_{2n}$ \cite{HR2005}.  

The \emph{twisted} versions of all the above monoids, $\PP_n^\tau$, $\PB_n^\tau$, $\B_n^\tau$, $\M_n^\tau$ and $\J_n^\tau$, are the corresponding submonoids of $\P_n^\tau$: thus, for example, the twisted Brauer monoid $\B_n^\tau$ has underlying set $\N\times\B_n$ and product~$\star$ given by \eqref{eq:star}.  In particular, the twisted Jones monoid $\J_n^\tau$ is known in the literature as the \emph{Kauffman monoid} and is denoted~$\K_n$ \cite{LF2006,BDP2002}.  Despite the above-mentioned isomorphism of $\PP_n$ and $\J_{2n}$, there is no such isomorphism between the twisted monoids $\PP_n^\tau$ and $\J_{2n}^\tau=\K_{2n}$.  Indeed, $\PP_n^\tau$ and $\K_{2n}$ do not have the same number of idempotents; see Tables \ref{tab-EJn-EKn} and \ref{tab-EPPn} in Section \ref{sect-values}.

The idempotents of the monoids $\P_n,\B_n,\PB_n$ (and their associated algebras and twisted versions) were classified and enumerated in \cite{DEEFHHL1}, and the purpose of the current article is to undertake the same program for their planar counterparts.
We conclude this subsection with a simple lemma.

\begin{lemma}\label{lem-rank}
For $\al\in\M_n$, the following are equivalent:
\begin{enumerate}[\rm (i)] \begin{multicols}{2}
\item $\al=\al^2$,
\item $\dom(\al)=\dom(\al^2)$,
\item $\codom(\al)=\codom(\al^2)$,
\item $\rank(\al)=\rank(\al^2)$.
\end{multicols}  
\end{enumerate}
\end{lemma}

\begin{proof}
(i)$\implies$(ii) and (i)$\implies$(iii).  These are obvious.

\pfitem{ii)$\implies$(iv) and (iii)$\implies$(iv}  If (ii) holds, then $\rank(\al)=|{\dom(\al)}| =|{\dom(\al^2)}|=\rank(\al^2)$, so that (iv) holds.  The other implication is dual.

\pfitem{iv)$\implies$(i}  Suppose $\rank(\al)=\rank(\al^2)$.  Every non-transversal of $\al$ is trivially a block of $\al^2$, so it remains to show that every transversal of $\al$ is a block of $\al^2$.  Let the transversals of $\al$ be $\{i_1,j_1'\},\ldots,\{i_r,j_r'\}$ where $r=\rank(\al)$ and $i_1<\cdots<i_r$, noting that this implies $\dom(\al)=\{i_1,\ldots,i_r\}$, $\codom(\al)=\{j_1,\ldots,j_r\}$ and $j_1<\cdots<j_r$ (the latter by planarity).  Since $\dom(\al^2)\sub\dom(\al)$ and $\codom(\al^2)\sub\codom(\al)$, every transversal of $\al^2$ is of the form $\{i_s,j_t'\}$ for some $s,t\in\{1,\ldots,r\}$.  Since $\rank(\al^2)=\rank(\al)=r$, there must be~$r$ such transversals of $\al^2$, and so these must be $\{i_1,j_{1\pi}'\},\ldots,\{i_r,j_{r\pi}'\}$ for some permutation $\pi$ of~$\{1,\ldots,r\}$.  Since $\al^2$ is planar, this permutation must be the identity, so it follows that $\al^2$ contains the transversals $\{i_1,j_1'\},\ldots,\{i_r,j_r'\}$, as required.
\end{proof}

\begin{rem}
We have not referred to Green's relations or the regular $*$-semigroup structure on $\M_n$, or any of the other diagram monoids we study, since neither plays a role in the theory developed in this section (though they do in the algorithms presented in Section \ref{sect-algorithms}).  Green's relations on $\M_n$, $\J_n$ and $\K_n$ are characterised in \cite[Theorem 2.4]{DEG2017}, \cite[Theorem 18]{Wilcox2007} and \cite[Theorem~5.1]{LF2006}, respectively, in terms of domains, ranks, and other parameters.  We will not need to know the exact formulations of these results, so we will not state them here, but it is worth noting that two elements of $\M_n$ are $\D$-related if and only if they have the same rank, and that $\M_n$ is $\H$-trivial.  These two facts lead to a simpler proof of Lemma \ref{lem-rank}, since for an element~$x$ of a finite semigroup, $x\D x^2\iff x\R x^2\iff x\L x^2\iff x\H x^2$ (see \cite[Theorems~A.2.4 and A.3.4]{RSbook}, for example); in particular, if the finite semigroup is~$\H$-trivial, then these are also equivalent to $x=x^2$.  The anti-involution~${}^*:\P_n\to\P_n$ that gives~$\P_n$, and hence all the submonoids considered in this article, a regular $*$-semigroup structure corresponds to reflecting (diagrams representing) elements of $\P_n$ in a horizontal axis midway between the two rows of vertices.
\end{rem}

\subsection{Interface graphs and characterisation of idempotents}\label{sect-interface}

A key role in our study is played by the so-called \emph{interface graph} of a Motzkin element.  In this subsection, we define these graphs, and show how they may be used to characterise the idempotents of $\M_n$, $\J_n$ and $\K_n$.

A block $A$ of a Motzkin element $\alpha\in\M_n$  is called an \emph{upper
hook} if $A\sub\bn$ and $|A|=2$, or an \emph{upper singleton} if $A\subseteq\bn$ and $|A|=1$.  \emph{Lower hooks} and \emph{lower singletons} are defined
analogously.

Let $\alpha\in\M_n$.  The \emph{interface graph} $\Gamma_\alpha$ is a vertex- and
edge-coloured graph defined as follows.  The vertex set of $\Gamma_\alpha$ is
simply $\bn$, and the colour $c(v)\in\Z_2\times\Z_2$ of a vertex $v\in\bn$ is
defined to be the column vector
\[
c(v)=\col ab \qquad\text{where}\qquad 
a=\begin{cases}
1 &\text{if $v\in\codom(\alpha)$}\\
0 &\text{otherwise}
\end{cases}
\AND
b=\begin{cases}
1 &\text{if $v\in\dom(\alpha)$}\\
0 &\text{otherwise.}
\end{cases}
\]
For each upper hook $\{i,j\}$ of $\alpha$, $\Gamma_\alpha$ has an edge $\{i,j\}$ coloured $-1$.  For each lower hook $\{k',l'\}$ of $\alpha$, $\Gamma_\alpha$ has an edge $\{k,l\}$ coloured $+1$.  (Note that $\Gamma_\alpha$ may have two edges between a pair of vertices, but these edges will always have opposite colours.)

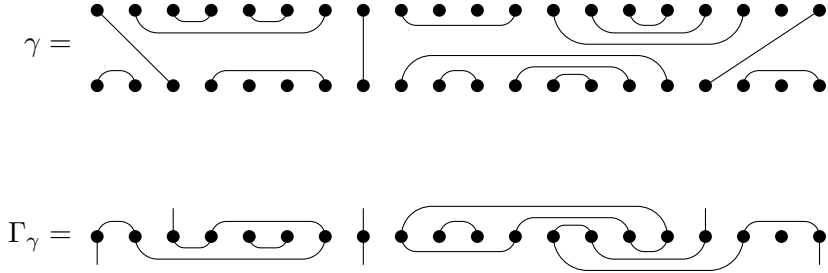
\begin{figure}[h]
\begin{center}
\begin{tikzpicture}[scale=.5]

\begin{scope}[shift={(0,0)}]	
\uvs{1,...,20}
\lvs{1,...,20}
\uarcx34{.3}
\uarcx56{.3}
\uarc9{12}
\uarc{15}{16}
\uarcx27{.6}
\uarcx{14}{17}{.6}
\uarcx{13}{18}{.9}
\darc12
\darc47
\darc{10}{11}
\darc{18}{20}
\darcx{13}{14}{.3}
\darcx{12}{15}{.5}
\darcx9{16}{.8}
\stline13
\stline88
\stline{20}{17}
\draw(0.6,1)node[left]{$\gamma=$};
\end{scope}

\begin{scope}[shift={(0,-4)}]	
\lvs{1,...,20}
\duarcx34{.3}
\duarcx56{.3}
\duarc9{12}
\duarc{15}{16}
\duarcx27{.6}
\duarcx{14}{17}{.6}
\duarcx{13}{18}{.9}
\darc12
\darc47
\darc{10}{11}
\darc{18}{20}
\darcx{13}{14}{.3}
\darcx{12}{15}{.5}
\darcx9{16}{.8}
\stlined1
\stlined8
\stlined{20}
\stlineu3
\stlineu8
\stlineu{17}
\draw(0.6,0)node[left]{$\Gamma_\gamma=$};
\end{scope}

\end{tikzpicture}
\end{center}
\vspace{-5mm}
\caption{A Motzkin element $\gamma\in\M_{20}$ (above) and its interface graph $\Gamma_\gamma$ (below).}
\label{fig:M20}
\end{figure}

When drawing the interface graph~$\Gamma_\alpha$ of a Motzkin element $\alpha\in\M_n$, we always draw the vertices in a horizontal row, in the order $1,\ldots,n$, increasing from left to right.  We draw the edges coloured~$+1$ or $-1$ above or below the line of vertices, respectively, as line segments:
\[
\begin{tikzpicture}[scale=.5]
\lvs{1,3,9,7}
\darc13
\duarc79
\draw(5,0)node {or};
\draw(9.5,-.2)node [right] {,};
\end{tikzpicture}
\]
with the ``height'' of such line segments chosen so that the diagram is planar (we can do this since~$\alpha\in\M_n$ is itself planar).  And we indicate the colour $c(v)=\col ab$ of a vertex $v\in\bn$ by drawing a small line segment above and/or below vertex $v$ if $a=1$ and/or $b=1$, respectively.  Roughly speaking, this picture of~$\Gamma_\alpha$ is obtained by cutting (a diagram representing) $\alpha$ in half, horizontally along the middle of the diagram, and then connecting the top half to the bottom half by identifying the two rows of vertices.
Figure~\ref{fig:M20} pictures a Motzkin element $\gamma\in\M_{20}$ and its interface graph~$\Gamma_\gamma$.

If $\alpha \in \M_n$, then every vertex in $\Gamma_{\alpha}$ has degree at most $2$. Hence, every connected component of~$\Gamma_{\alpha}$ is either a cycle or a path; we regard a singleton component of $\Gamma_\alpha$ as a path of length $0$.  A vertex cannot be the endpoint of two edges with the same colour, because blocks of $\alpha$ have size at most $2$, so it follows that the edges along any path in $\Gamma_\alpha$ alternate in colour; in particular, all cycles have even length.  It is also apparent that a vertex of degree $2$ can only be coloured $\col00$.  Now consider a path component $v_1-v_2-\cdots-v_k$ of $\Gamma_\alpha$.  As above, $c(v_2)=\cdots=c(v_{k-1})=\col00$.  We call the path \emph{inactive} if also $c(v_1)=c(v_k)=\col00$.  We call the path \emph{active} if either $k=1$ and $c(v_1)=\col11$ or else $k\geq2$ and $c(v_1),c(v_k)\in\left\{\col01,\col10\right\}$.  Otherwise, we say the path is \emph{mixed}.

The next result characterises the idempotents of $\M_n$ in terms of interface graphs.  Of crucial importance is the fact that the product graph $\Pi(\alpha,\alpha)$ contains an isomorphic copy of $\Gamma_\alpha$ in the middle layer.  If $C$ is a connected component of a graph $\Gamma$, so that $C$ is itself a graph, we will slightly abuse notation and identify~$C$ with its set of vertices, so we sometimes write ``$u\in C$'' to mean ``$u$ is a vertex of $C$''.

\begin{prop}\label{prop-idempotents}
A Motzkin element $\alpha \in \M_n$ is an idempotent if and only if every connected component of the interface graph $\Gamma_{\alpha}$ is one of: 
  \begin{enumerate}[\rm (i)]\begin{multicols}{3}
    \item 
      a cycle,
    \item 
      an inactive path, or
    \item 
      an active path of even length.
\end{multicols}  \end{enumerate}
\end{prop}

\begin{proof}
($\Leftarrow$).  Suppose first that all components of $\Gamma_\alpha$ are of types (i)--(iii).   By Lemma \ref{lem-rank}, and since clearly $\rank(\al^2)\leq\rank(\al)$, to show that $\alpha^2=\alpha$, it is enough to show that $\rank(\alpha^2)\geq\rank(\alpha)$.  

Suppose the components of $\Gamma_\alpha$ of type (iii) are $C_1,\ldots,C_r$.  For each $i$, let the endpoints of $C_i$ be $u_i$ and~$v_i$.  Since $C_i$ is of even length, we may assume that the bottom coordinate of $c(u_i)$ and the top coordinate of~$c(v_i)$ are both $1$, even if $u_i=v_i$.  

As noted above, the colour of any vertex from a component of types (i) or (ii) is $\col00$; this is also the case for any interior vertices of the components of type (iii).  It follows that $\dom(\alpha)=\{u_1,\ldots,u_r\}$ and $\codom(\alpha)=\{v_1,\ldots,v_r\}$, so that $\rank(\alpha)=r$.  It also follows that there is a permutation $\pi$ of $\{1,\ldots,r\}$ such that $\{u_1,v_{\pi(1)}'\},\ldots,\{u_r,v_{\pi(r)}'\}$ are the transversals of $\alpha$.

Now fix some $1\leq i\leq r$.  Since $C_i$ is a path from $u_i$ to $v_i$, it follows that the product graph $\Pi(\alpha,\alpha)$ has a path from $u_i''$ to $v_i''$.  Since $\{u_{\pi^{-1}(i)},v_i'\}$ and $\{u_i,v_{\pi(i)}'\}$ are transversals of $\alpha$, $\Pi(\alpha,\alpha)$ contains the edges~$\{u_{\pi^{-1}(i)},v_i''\}$ and $\{u_i'',v_{\pi(i)}'\}$.  So $u_{\pi^{-1}(i)}$ and $v_{\pi(i)}'$ are connected by a path in $\Pi(\alpha,\alpha)$, and it follows that~$\{u_{\pi^{-1}(i)},v_{\pi(i)}'\}$ is a transversal of $\alpha^2$.  Thus, $\dom(\al^2)\supseteq\{u_{\pi^{-1}(1)},\ldots,u_{\pi^{-1}(r)}\}=\{u_1,\ldots,u_r\}$, and it follows that ${\rank(\alpha^2)\geq r=\rank(\alpha)}$, as required.

\pfitem{$\Rightarrow$}  For this implication, we prove the contrapositive.  Suppose $\Gamma_\alpha$ contains a component not of types~${\text{(i)--(iii)}}$.  Then this component must be either
  \begin{enumerate}\begin{multicols}{2}
    \item[(iv)]
      an active path of odd length, or
    \item[(v)]
      a mixed path.
\end{multicols}  \end{enumerate}
Suppose first that $C$ is a component of $\Gamma_\alpha$ of type (iv), and let $u,v$ be the endpoints of $C$.  Since $C$ is of odd length, it follows that $u\not=v$, and that $c(u)=c(v)=\col10$ or $\col01$.  We suppose the latter is the case (the proof for the former is similar).  So $u,v\in\dom(\alpha)$ belong to (distinct) transversals of $\al$: $\{u,x'\}$ and $\{v,y'\}$, say.  So $\{u'',x'\}$ and $\{v'',y'\}$ are edges in the product graph $\Pi(\alpha,\alpha)$.  Since $C$ gives rise to a path from $u''$ to $v''$ in $\Pi(\alpha,\alpha)$, it follows that $\{x',y'\}$ is a block of $\alpha^2$.  But $\{x',y'\}$ is not a block of $\alpha$, since $x,y\in\codom(\alpha)$, so it follows that $\alpha^2\not=\alpha$.  

On the other hand, suppose $C$ is a component of $\Gamma_\alpha$ of type (v), and let $u,v$ be the endpoints of $C$.  Suppose $c(v)=\col00$, so that $c(u)=\col10$ or $\col01$.  Again, we just consider the latter case.  So $u$ belongs to a transversal $\{u,x'\}$ of $\alpha$.  Since $c(v)=\col00$, it follows that the connected component in the product graph $\Pi(\alpha,\alpha)$ containing $x'$ is $\{x'\}\cup\set{y''}{y\in C}$.  We deduce that $\{x'\}$ is a block of $\alpha^2$, and it again follows that $\alpha^2\not=\alpha$, since $x\in\codom(\alpha)$.  This completes the proof.
\end{proof}

The interface graph $\Gamma_\alpha$ of a Jones element $\alpha\in\J_n$ can only have cycles and active paths, since all blocks of $\alpha$ are of size $2$.  So we may immediately deduce from Proposition \ref{prop-idempotents} the following characterisation of Jones idempotents.

\begin{cor}\label{cor-jones-idempotents}
A Jones element $\alpha \in \J_n$ is an idempotent if and only if every connected component of $\Gamma_{\alpha}$ is a cycle or an active path of even length. \qed
\end{cor}

The characterisation of idempotents in the twisted Motzkin and Jones monoids, $\M_n^\tau$ and $\J_n^\tau=\K_n$, is as follows.

\begin{cor}\label{cor-kauffman-idempotents}
A twisted Motzkin element $(i,\alpha) \in \M_n^\tau$ is an idempotent if and only if $i=0$ and every connected component of $\Gamma_{\alpha}$ is an active path of even length.  Consequently, ${E(\M_n^\tau)=E(\K_n)}$.
\end{cor}

\begin{proof}
Note that $(i,\alpha)\star(i,\alpha)=(2i+\tau(\alpha,\alpha),\alpha^2)$, so $(i,\alpha)\in E(\M_n^\tau)$ if and only if $i=0$, $\tau(\alpha,\alpha)=0$ and $\alpha\in E(\M_n)$.  Cycles and inactive paths in the interface graph $\Gamma_\alpha$ correspond to floating components in the product graph $\Pi(\alpha,\alpha)$, so it follows from Proposition \ref{prop-idempotents} that [$\tau(\alpha,\alpha)=0$ and $\alpha\in E(\M_n)$] is equivalent to~$\Gamma_\alpha$ having only active paths of even length.  The last assertion of the lemma follows quickly.
\end{proof}

\begin{rem}\label{rem:TL1}
Corollary \ref{cor-kauffman-idempotents} also applies to the \emph{Temperley-Lieb} and \emph{Motzkin algebras}.  If $S$ is any of the monoids $\P_n$, $\PB_n$, $\B_n$, $\PP_n$, $\M_n$ or $\J_n$, and if $\F$ is a field with some fixed element $\xi\in\F$, then we may form the \emph{twisted semigroup algebra} $\F^\xi[S]$, as in \cite{Wilcox2007}; these are the partition \cite{Martin1994,Jones1994_2}, partial Brauer \cite{MarMaz2014}, Brauer \cite{Brauer1937}, planar partition \cite{Jones1994_2}, Motzkin \cite{BH2014} and Temperley-Lieb \cite{TL1971} algebras, respectively.  The algebra~$\F^\xi[S]$ has basis~$S$, and multiplication $\circ$ defined on basis elements $\al,\be\in S$ (and extended linearly) by $\al\circ\be=\xi^{\tau(\al,\be)}(\al\be)$.  If $\xi$ is not a root of unity or if it is an $M$th root of unity where $M>n$ (the so-called \emph{generic case}), then an element $\al\in S$ satisfies $\al=\al\circ\al$ if and only if $\tau(\al,\al)=0$ and $\al=\al^2$ in $S$; thus, in this case, Corollary~\ref{cor-kauffman-idempotents} shows that an element $\al\in\M_n$ is an idempotent basis element of the Motzkin algebra if and only if $i=0$ and every connected component of $\Gamma_{\alpha}$ is an active path of even length.  
As in \cite[Section~6]{DEEFHHL1}, if $\xi$ is an $M$th root of unity with $M\leq n$, then a Motzkin element $\al\in\M_n$ is an idempotent basis element of $\F^\xi[\M_n]$ if and only every connected component of $\Ga_\al$ is one of types (i)--(iii) as listed in Proposition \ref{prop-idempotents}, with the combined number of components of types (i)--(ii) being a multiple of~$M$.
%
On the other hand, if $\al$ is an arbitrary idempotent of~$S$ (where $S$ is any of the above diagram monoids), and if $\xi\not=0$, then $\xi^{-m(\al,\al)}\al$ is an idempotent of $\F^\xi[S]$; we thank Zajj Daugherty for this last observation.
\end{rem}

\subsection{A mapping on \boldmath{$E(\M_n)$} and an enumeration method}\label{sect-mapping}

Now that we have characterised the idempotents of $\M_n$, $\J_n$ and $\K_n$, we wish to enumerate them.  In this subsection, we describe a method for doing so.  We make crucial use of a map $D:E(\M_n)\to E(\M_n)$ to be defined shortly, and the interface graphs defined in Subsection \ref{sect-interface}.  This map also played a crucial role in the classification of \emph{congruences} on $\J_n$ and $\B_n$ in \cite{EMRT2018}; a different, but closely-related, map was used for~$\M_n$,~$\P_n$,~$\PB_n$ and~$\PP_n$.

\begin{lemma}\label{lem-paths}
Suppose $\alpha\in E(\M_n)$ and $\{i,j'\}$ is a transversal of $\alpha$.  Then $\Gamma_\alpha$ contains an active path of even length with $i$ and $j$ as its endpoints.  In particular, if $\dom(\alpha)=\{i_1,\ldots,i_r\}$ where $r=\rank(\alpha)$, then $i_1,\ldots,i_r$ belong to $r$ distinct connected components of $\Gamma_\alpha$.
\end{lemma}

\begin{proof}
Since $\alpha=\alpha^2$, $\{i,j'\}$ is a transversal of $\al^2$, so there is a path from $i$ to $j'$ in the product graph~$\Pi(\alpha,\alpha)$.  Since $\{i,j''\}$ and $\{i'',j'\}$ are both edges of $\Pi(\alpha,\alpha)$, it follows that there is a path from $j''$ to $i''$ in~$\Pi(\alpha,\alpha)$ involving only vertices in the middle row; this gives rise to a path from $i$ to $j$ in $\Gamma_\alpha$.  Since $i\in\dom(\alpha)$ and $j\in\codom(\alpha)$, it follows that $c(i)\not=\col00$ and $c(j)\not=\col00$, so that $i$ and $j$ are indeed the endpoints of this path.  Since this path is active, Proposition \ref{prop-idempotents} tells us that it is of even length.  This proves the first assertion of the lemma.

For the second assertion, suppose the transversals of $\al$ are $\{i_1,j_1'\},\ldots,\{i_r,j_r'\}$, and let $1\leq k<l\leq r$.  We must show that $i_k$ and $i_l$ belong to different connected components of $\Ga_\al$.  To do so, suppose to the contrary that $i_k$ and $i_l$ belong to the same component.  By the previous paragraph, the component of $\Ga_\al$ containing $i_k$ is a path from $i_k$ to $j_k$, and the component containing $i_l$ is a path from $i_l$ to $j_l$.  Thus, since we have assumed these are the same components, and since $i_k\not=i_l$, we must have $i_k=j_l$ and $i_l=j_k$.  But then $\{i_k,i_l'\}$ and $\{i_l,i_k'\}$ are both transversals of $\al$, contradicting planarity.
\end{proof}

To define the mapping $D:E(\M_n)\to E(\M_n)$, we first define a map $d:E(\M_n)\to\M_n$.  
Before we do this, note first that the set of transversals of a Motzkin element $\alpha\in\M_n$ inherits an obvious total ordering from the natural ordering on $\bn$.  For example, the transversals of the Motzkin element $\gamma\in\M_{20}$ pictured in Figure \ref{fig:M20} are ordered by $\{1,3'\}<\{8,8'\}<\{20,17'\}$.  So we may speak of the first and second transversals of~$\alpha$, and so on.

\begin{figure}[h]
\begin{center}
\begin{tikzpicture}[scale=.5]

\begin{scope}[shift={(0,0)}]	
\bluebox1{1.7}42
\bluebox6{1.7}82
\bluebox103{.3}
\bluebox509{.3}
\bluetrap{11}0{14}0{14}2{10}2
\uvs{1,4,5,6,8,9,10,14}
\lvs{1,3,4,5,9,10,11,14}
\stlines{5/4,9/10}

\stlinedots{10/11,14/14}
\uvertlabels{1/1,5/i,9/k,14/n}
\lvertlabels{1/1,4/j,10/l,14/n}
\udotteds{1/4,6/8,10/14}
\ldotteds{1/3,5/9,11/14}
\draw(.5,1)node[left]{$\alpha=$};
\end{scope}

\begin{scope}[shift={(0,-5)}]	
\bluebox1{1.7}42
\bluebox6{1.7}82
\bluebox103{.3}
\bluebox509{.3}
\bluetrap{11}0{14}0{14}2{10}2
\uvs{1,4,5,6,8,9,10,14}
\lvs{1,3,4,5,9,10,11,14}

\uarcx59{.7}
\darcx4{10}{.7}

\stlinedots{10/11,14/14}
\uvertlabels{1/1,5/i,9/k,14/n}
\lvertlabels{1/1,4/j,10/l,14/n}
\udotteds{1/4,6/8,10/14}
\ldotteds{1/3,5/9,11/14}
\draw(.5,1)node[left]{$d(\alpha)=$};
\end{scope}

\end{tikzpicture}
\end{center}
\vspace{-5mm}
\caption{The map $d:E(\M_n)\to E(\M_n)$.  Shaded regions of $\alpha$ are assumed to be identical to the corresponding shaded regions of $d(\alpha)$.}
\label{fig:d}
\end{figure}
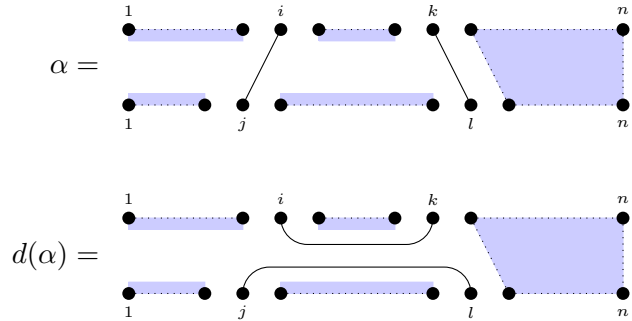

Now let $\alpha\in E(\M_n)$.  If $\rank(\alpha)\leq1$, then we define $d(\alpha)=\alpha$.  Otherwise, if $\{i,j'\}$ and $\{k,l'\}$ are the first two transversals of $\alpha$, then we define $d(\alpha)$ to be the element of $\M_n$ obtained from $\alpha$ by replacing these two transversals by the upper and lower hooks $\{i,k\}$ and $\{j',l'\}$.  Note that blocks of $d(\alpha)$ trivially have size~$\leq2$, while the planarity of $\alpha$ ensures that of $d(\alpha)$; see Figure~\ref{fig:d}.  Although it is not readily apparent that~$d(\alpha)$ is necessarily an idempotent, we will soon see that it is.

The next result gathers some important properties of the $d$ map, including the fact that $d$ does indeed map into $E(\M_n)$.  
Let $A$ be an upper non-transversal of a Motzkin element $\alpha\in\M_n$.  We say that $A$ is \emph{to the left} of a transversal $\{i,j'\}$ of $\alpha$ if $\max(A)<i$.  Similarly, we say a lower non-transversal $B'$ is to the left of $\{i,j'\}$ if $\max(B)<j$.  We say that $A$ is \emph{nested} (in $\alpha$) if there exists an upper hook $\{k,l\}$ of $\alpha$ such that~$k<\min(A)$ and $\max(A)<l$; otherwise, we say that $A$ is \emph{unnested}.  We define nested and unnested lower non-transversals analogously.  An \emph{outer hook} of $\alpha$ is defined to be an upper or lower unnested hook of $\alpha$ that is to the left of any transversal of $\alpha$.  For the statements of parts (ii) and (iii) of the next lemma, recall that we informally identify a connected component of a graph with its underlying vertex set.

\begin{lemma}\label{lem-d}
Let $\alpha\in E(\M_n)$ with $\rank(\alpha)\geq2$, and suppose $\{i,j'\}$ and $\{k,l'\}$ are the first two transversals of $\alpha$.  Suppose $C_1,C_2,C_3,\ldots,C_s$ are the connected components of $\Gamma_\alpha$, where $i\in C_1$ and $k\in C_2$.  Then
\begin{enumerate}[\rm (i)]
\item $\{i,k\}$ and $\{j',l'\}$ are outer hooks of $d(\alpha)$,
\item $C_1\cup C_2,C_3,\ldots,C_s$ are the connected components of $\Gamma_{d(\alpha)}$,
\item $i,j,k,l$ belong to $C_1\cup C_2$, and this is a cycle component of $\Gamma_{d(\alpha)}$, and
\item $d(\alpha)\in E(\M_n)$.
\end{enumerate}
\end{lemma}

\begin{proof}
(i).  Planarity of $\alpha$, and the fact that there are no other transversals between $\{i,j'\}$ and $\{k,l'\}$ ensures that $\{i,k\}$ and $\{j',l'\}$ are unnested in $d(\alpha)$; cf.~Figure \ref{fig:d}.  The fact that $\{i,j'\}$ and $\{k,l'\}$ are the \emph{first} two transversals of $\alpha$ ensures that $\{i,k\}$ and $\{j',l'\}$ are to the left of any transversal of $d(\alpha)$.  

\pfitem{ii) and (iii}  By Lemma \ref{lem-paths}, $C_1$ is a path in $\Gamma_\alpha$ from $i$ to $j$, and $C_2$ a path from $k$ to $l$.  The only change from $\Ga_\al$ to $\Ga_{d(\al)}$ is the addition of the edges $\{i,k\}$ and $\{j,l\}$, coloured $-1$ and $+1$, respectively, and the recolouring of the vertices $i,j,k,l$ (four $1$'s are changed to $0$'s, regardless of whether these are four distinct vertices).  Since the addition of these edges joins $C_1$ and $C_2$ into a single cycle component of $\Ga_{d(\al)}$, (iii) follows.  Since no other components of $\Gamma_\alpha$ are modified,~(ii) also follows.


\pfitem{iv} Since $\alpha\in E(\M_n)$, the components $C_3,\ldots,C_s$ are all of the forms specified in Proposition~\ref{prop-idempotents}.  Since $C_1\cup C_2$ is a cycle, (iv) follows.
\end{proof}

\begin{figure}[h]
\begin{center}
\begin{tikzpicture}[scale=.5]

\begin{scope}[shift={(0,0)}]	
\uarcx34{.3}
\uarcx56{.3}
\uarc9{12}
\uarc{15}{16}
\uarcx27{.6}
\uarcx{14}{17}{.6}
\uarcx{13}{18}{.9}
\darc12
\darc47
\darc{10}{11}
\darc{18}{20}
\darcx{13}{14}{.3}
\darcx{12}{15}{.5}
\darcx9{16}{.8}
\uarcxred181
\darcxred38{.7}

\stline{20}{17}
\draw(0.6,1)node[left]{$\delta=d(\gamma)=$};
\uvs{1,...,20}
\lvs{1,...,20}
\end{scope}

\begin{scope}[shift={(0,-4)}]	
\duarcx34{.3}
\duarcx56{.3}
\duarc9{12}
\duarc{15}{16}
\duarcx27{.6}
\duarcx{14}{17}{.6}
\duarcx{13}{18}{.9}
\darc12
\darc47
\darc{10}{11}
\darc{18}{20}
\darcx{13}{14}{.3}
\darcx{12}{15}{.5}
\darcx9{16}{.8}
\duarcxred181
\darcxred38{.7}

\stlined{20}

\stlineu{17}
\draw(0.6,0)node[left]{$\Gamma_\delta=$};
\lvs{1,...,20}
\end{scope}

\end{tikzpicture}
\end{center}
\vspace{-5mm}
\caption{The Motzkin element $\delta=d(\gamma)\in\M_{20}$ (above) and its interface
graph $\Gamma_\delta$ (below), where $\gamma\in\M_{20}$ is pictured in Figure
\ref{fig:M20}.  New edges are drawn in red.}
\label{fig:M20cont}
\end{figure}
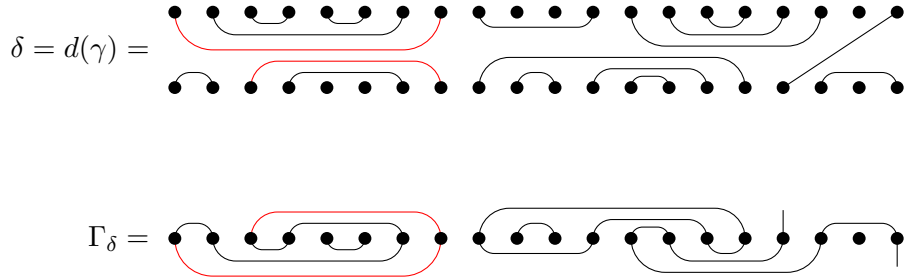

For any $\alpha\in E(\M_n)$, the sequence $\alpha,d(\alpha),d^2(\alpha),\ldots$ eventually terminates in an idempotent of rank $0$ or~$1$, depending on the parity of $\rank(\alpha)$, and we write $D(\alpha)$ for this idempotent.  In fact, $D(\alpha)=d^s(\alpha)$, where $s=\lfloor \rank(\alpha)/2\rfloor$; we consider $d^0$ to be the identity mapping.
For example, consider the Motzkin element $\gamma\in\M_{20}$ pictured in Figure \ref{fig:M20}.  Here we have $\rank(\gamma)=3$, so that $D(\gamma)=d(\gamma)$; we have pictured $\delta=d(\gamma)$ and $\Gamma_\delta$ in Figure \ref{fig:M20cont}.

From now on, we will write
\begin{align*}
\De(\M_n) = D(E(\M_n)) = \set{D(\alpha)}{\alpha\in E(\M_n)} &\AND
\De(\J_n) = D(E(\J_n)) = \set{D(\alpha)}{\alpha\in E(\J_n)}.
\intertext{Since $D$ acts as the identity on idempotents of rank $\leq1$, and since all Jones elements of minimal rank are idempotents, it follows that}
 \De(\M_n) = \set{\alpha\in E(\M_n)}{\rank(\alpha)\leq1} &\AND  \De(\J_n) = \set{\alpha\in \J_n}{\rank(\alpha)\leq1}.
\end{align*}
Note that there are elements in $\J_n$ of rank $0$ or $1$, depending on the parity of $n$,
but not both. 

The next result shows that enumeration of $E(\M_n)$ and $E(\J_n)$ reduces to
the enumeration of preimages under the~$D$ map.  This latter task is itself
quite difficult; Section \ref{sect-algorithms} and the remainder of Section \ref{sect-theory} are devoted to achieving it.

\begin{lemma}\label{lem-EMJ}
If $S$ is one of $\M_n$ or $\J_n$, then $\displaystyle{|E(S)| = \sum_{\alpha\in \De(S)}|D^{-1}(\alpha)|}$.

\end{lemma}

\begin{proof}
  The statement about $|E(\M_n)|$ is clear.  If $\alpha\in E(\M_n)$, then
  Corollary~\ref{cor-jones-idempotents} and Lemma~\ref{lem-d} imply that
  $\alpha\in\J_n$ if and only if $D(\alpha)\in\J_n$.  The statement about $|E(\J_n)|$
  follows.
\end{proof}

For an arbitrary Motzkin element $\alpha\in\M_n$, we write $\Theta(\alpha)$ for the set of all cycle components of the interface graph $\Gamma_\alpha$ containing at least one edge corresponding to an upper outer hook of $\alpha$ and at least one edge corresponding to a lower outer hook of $\alpha$; outer hooks were defined before Lemma \ref{lem-d}.  For $\theta\in\Theta(\alpha)$, let
\begin{align*}
U_\theta(\alpha) &= \bigset{\{i,j\}}{\text{$\{i,j\}$ is an outer hook of $\alpha$ and $i,j\in\theta$}},\\
L_\theta(\alpha) &= \bigset{\{i',j'\}}{\text{$\{i',j'\}$ is a outer hook of $\alpha$ and $i,j\in\theta$}},
\end{align*}
and write $u_\theta(\alpha)=|U_\theta(\alpha)|$ and $l_\theta(\alpha)=|L_\theta(\alpha)|$.  
For example, we have $\Theta(\gamma)=\emptyset$ for $\gamma\in\M_{20}$ from Figure~\ref{fig:M20} (note that $\gamma$ has only one outer hook).  However, with $\delta=d(\gamma)\in\M_{20}$ from Figure~\ref{fig:M20cont}, we have $\Theta(\delta)=\{\theta_1,\theta_2\}$, where $\theta_1=\{1,2,3,4,7,8\}$ and $\theta_2=\{9,12,15,16\}$, and 
\[
U_{\theta_1}(\delta) = \big\{ \{1,8\} \big\} , \ \ \  \
L_{\theta_1}(\delta) = \big\{ \{1',2'\},\{3',8'\}\big\} , \ \ \ \
U_{\theta_2}(\delta) = \big\{ \{9,12\} \big\} , \ \  \ \
L_{\theta_2}(\delta) = \big\{ \{9',16'\}\big\} .
\]
So $l_{\theta_1}(\delta)=2$, while $u_{\theta_1}(\delta)=u_{\theta_2}(\delta)=l_{\theta_2}(\delta)=1$.
The next result shows why the sets we have just defined are important.

\begin{lemma}\label{lem-d-1}
Suppose $\alpha\in E(\M_n)$ is such that $\Theta(\alpha)\not=\emptyset$.  Let $\theta\in\Theta(\alpha)$, and let $\{i,k\}\in U_\theta(\alpha)$ and $\{j',l'\}\in L_\theta(\alpha)$.  Let $\beta$ be obtained from $\alpha$ by replacing the blocks $\{i,k\}$ and $\{j',l'\}$ by $\{i,j'\}$ and $\{k,l'\}$.  Then $\beta\in E(\M_n)$ and $d(\beta)=\alpha$.  Further, every element of $d^{-1}(\alpha)\setminus\{\alpha\}$ may be constructed in this way, for some $\theta\in\Theta(\alpha)$ and some pair of edges from $U_\theta(\alpha)\times L_\theta(\alpha)$.
\end{lemma}

\begin{proof}
Since $\{i,k\}$ and $\{j',l'\}$ are outer hooks of $\alpha$, it follows that $\beta\in\M_n$.  To show that $\beta$ is an idempotent, we need to check that each component of $\Gamma_\beta$ is of one of the forms specified in Proposition~\ref{prop-idempotents}.  By construction, the components of $\Gamma_\alpha$ other than $\theta$ are still components of $\Gamma_\beta$, and these must all be of the specified form, since $\alpha\in E(\M_n)$.  But $\theta$, a cycle component of $\Gamma_\alpha$, is split into two active path components of $\Gamma_\beta$.  To complete the proof that $\beta$ is an idempotent, it remains to show that these paths are of even length.  But this follows quickly from the fact that we are removing an edge coloured $+1$ and an edge coloured $-1$ from a cycle whose edge colours alternate between $+1$ and $-1$.

Since $\{i,k\}$ and $\{j',l'\}$ are to the left of any transversals of $\alpha$, as they are outer hooks of $\alpha$, it follows that $\{i,j'\}$ and $\{k,l'\}$ are the first two transversals of $\beta$, and then it follows immediately that $d(\beta)=\alpha$.

Finally, suppose $\gamma\in E(\M_n)$ is such that $d(\gamma)=\alpha$, and let the first two transversals of $\gamma$ be $\{u,v'\}$ and $\{x,y'\}$, respectively.  By Lemma \ref{lem-d}, $\{u,x\}$ and $\{v',y'\}$ are outer hooks of $d(\gamma)=\alpha$, and $u,v,x,y$ belong to the same cycle component of $\Gamma_{d(\gamma)}=\Gamma_\alpha$.  If we denote this cycle component by $\sigma$, then $\{u,x\}\in U_\sigma(\alpha)$ and $\{v',y'\}\in L_\sigma(\alpha)$, and we see that $\gamma$ is constructed in the manner described in the lemma, with respect to $\sigma$, $\{u,x\}$ and $\{v',y'\}$.
\end{proof}

Lemma \ref{lem-d-1} gives information about preimages under the $d$ map.  In order to extend this to preimages under the $D$ map, we require the next two intermediate lemmas.  The first, Lemma~\ref{lem-geometry}, concerns curves in the plane, and the second, Lemma~\ref{lem-upper}, applies this to the situation in which the curves are part of the interface graph of a Motzkin element.

\begin{lemma}\label{lem-geometry}
Let $A,B,C,D$ be distinct points on the $x$-axis, with $x$-coordinates $a<b<c<d$, respectively.  Suppose $\C_1$ and $\C_2$ are smooth non-self-intersecting curves in the plane such that
\begin{enumerate}[\rm (i)]
\item $\C_1$ joins $A$ to $C$, while $\C_2$ joins $B$ to $D$,
\item apart from the endpoints stated above, both curves are contained in the region $a<x<d$, and
\item $\C_1$ and $\C_2$ never go below the points $B$ or $C$: that is, no point $(x,y)$ on either curve satisfies $[x=b$ and $y<0]$ or $[x=c$ and $y<0]$.
\end{enumerate}
Then $\C_1$ and $\C_2$ intersect.
\end{lemma}

\begin{proof}
Consider the curve $\C$ obtained from $\C_1$ by adding the positive half $\L_1$ of the line $x=a$ and the negative half $\L_2$ of the line $x=c$, as shown in Figure \ref{fig:C1}.  By the stated assumptions, $\C$ has no self-intersections, and so divides the plane into two regions: one containing $B$ and one containing $D$.  But $\C_2$ joins $B$ to $D$, so it follows that that $\C_2$ and $\C$ intersect.  
Assumptions (ii) and (iii), respectively, tell us that~$\C_2$ does not intersect $\L_1$ or $\L_2$.  So $\C_2$ must intersect $\C_1$.
\end{proof}

\begin{figure}[h]
\begin{center}
\newcommand{\scale}{1.5}
\begin{tikzpicture}[scale=1]
\draw[gray, <->] (\scale*0,0)--(\scale*10,0);
\node () at (\scale*2,-.1) [below] {$A$};
\node () at (\scale*4,-.1) [below] {$B$};
\node () at (\scale*6,-.1) [below left] {$C$};
\node () at (\scale*8,-.1) [below] {$D$};
\node[blue] () at (\scale*2,3.8) [left] {$\L_1$};
\node[blue] () at (\scale*6,-2) [right] {$\L_2$};
\node[red] () at (\scale*4,3.5) [right] {$\C_1$};
\draw[blue,ultra thick, ->] (\scale*2,0)--(\scale*2,4);
\draw[blue, ultra thick, ->] (\scale*6,0)--(\scale*6,-2);
\draw[red, ultra thick] plot [smooth] coordinates {(\scale*2,0) (\scale*2.3,.3) (\scale*2.7,-.5) (\scale*3.2,.4) (\scale*2.9,1.2) (\scale*4.3,1.8) (\scale*3.7,2.8) (\scale*4.2,3.3) (\scale*5,1.6) (\scale*3.7,.5) (\scale*5.5,-.3) (\scale*6,2.3) (\scale*7.5,.7) (\scale*6.8,-.5) (\scale*6.3,.5) (\scale*5.8,.2) (\scale*6,0)};
\foreach \x in {\scale*2,\scale*4,\scale*6,\scale*8} {\fill (\x,0)circle(.1);}
\end{tikzpicture}
\end{center}
\vspace{-8mm}
\caption{The curve ${\C}={\L_1}\cup{\C_1}\cup{\L_2}$ from the proof of Lemma \ref{lem-geometry}.}
\label{fig:C1}
\end{figure}

Recall that we are identifying a connected component $C$ of a graph on vertex set $\bn$ with the underlying vertex set of $C$.  In this way, we may also write $\min(C)$ to mean the vertex of $C$ with minimum value in the natural ordering on $\bn$.

\newpage

\begin{lemma}\label{lem-upper}
Suppose $\alpha\in E(\M_n)$ and $\theta_1,\theta_2\in\Theta(\alpha)$, where $\min(\theta_1)<\min(\theta_2)$.  
\begin{enumerate}[\rm (i)]
\item If $\{i,j\}\in U_{\theta_1}(\alpha)$ and $\{k,l\}\in U_{\theta_2}(\alpha)$ with $i<j$ and $k<l$, then $j<k$.
\item If $\{i',j'\}\in L_{\theta_1}(\alpha)$ and $\{k',l'\}\in L_{\theta_2}(\alpha)$ with $i<j$ and $k<l$, then $j<k$.
\end{enumerate}
\end{lemma}

\begin{proof}
We just prove (i), as (ii) is dual.
Suppose to the contrary that $j>k$.  Since $\{i,j\}$ and $\{k,l\}$ are unnested blocks of $\alpha$, and since $\alpha$ is planar, it follows that 
$k<l<i<j$.  In what follows, we consider the cycles $\theta_1$ and $\theta_2$ as (closed, non-self-intersecting) curves in the plane, with each vertex $v\in\bn$ drawn at the point $(v,0)$, and with edges drawn in the usual way to join vertices as:
\[
\begin{tikzpicture}[scale=.5]
\lvs{1,3,9,7}
\darc13
\duarc79
\draw(5,0)node {or};
\draw(9.5,-.2)node [right] {.};
\end{tikzpicture}
\]
We first claim that $\max(\theta_2)<\max(\theta_1)$.  Indeed, suppose to the contrary that $\max(\theta_2)>\max(\theta_1)$.  Put $a=\min(\theta_1)$, $b=l$, $c=i$ and $d=\max(\theta_2)$, and let $A=(a,0)$, $B=(b,0)$, and so on.  Note that $\theta_1$ is the union of two paths joining $A$ and $C$; let $\C_1$ be either of these paths.  Similarly, $\theta_2$ is the union of two paths joining $B$ and $D$; let $\C_2$ be either of these paths.  It is easy to check that conditions (i)--(iii) of Lemma~\ref{lem-geometry} are satisfied, using the fact that $\{i,j\}$ and $\{k,l\}$ are unnested to verify condition (iii).  It follows that $\C_1$ and $\C_2$, and hence $\theta_1$ and $\theta_2$, intersect, a contradiction.  This completes the proof of the claim that $\max(\theta_2)<\max(\theta_1)$.  It follows that $\min(\theta_1)<\min(\theta_2)<\max(\theta_2)<\max(\theta_1)$. 

Now, $\th_1$ is also the union of two paths joining $(\min(\th_1),0)$ and $(\max(\th_1),0)$; let $\C_3$ be either of these paths.
So $\C_3$ is a non-self-intersecting curve in the plane
and, apart from its endpoints, it lies in the region $\min(\theta_1)<x<\max(\theta_1)$.  In particular, it divides the region $\min(\theta_1)<x<\max(\theta_1)$ into upper and lower regions.  Since $\theta_1$ and $\theta_2$ do not intersect, $\theta_2$ is contained wholly within one of these two regions.  Since $\{k,l\}$ is unnested, it lies in the lower region and, hence, it follows that $\theta_2$ is contained in this lower region.  But then every edge of $\theta_2$ lies under the curve $\C_3\subseteq\theta_1$.  It follows that every lower hook of $\al$ corresponding to an edge of $\th_2$ is nested in $\al$, so that $L_{\theta_2}(\alpha)=\emptyset$, contradicting the assumption that $\theta_2\in\Theta(\alpha)$.
\end{proof}

\begin{rem}
As the last paragraph of the above proof indicates, the assumption that $\theta_2$ has outer upper \emph{and} lower hooks is necessary to prove the conclusion of Lemma \ref{lem-upper}(i).  Indeed, consider the Jones idempotent 
$\al=\big\{ \{1,2\}, \{3,4\}, \{5,6\}, \{1',6'\}, \{2',5'\}, \{3',4'\} \big\}\in E(\J_6)$.  
The interface graph $\Ga_\al$ has two connected components: $\th_1=\{1,2,5,6\}$ and $\th_2=\{3,4\}$.  Both are cycles, and $\min(\th_1)<\min(\th_2)$, yet $\{5,6\}$ and $\{3,4\}$ are upper outer hooks of $\th_1$ and $\th_2$, respectively, and we do not have $6<3$.  However, while $\th_1$ does belong to $\Th(\al)$, $\th_2$ does not.
\end{rem}

We are now ready to combine the preceeding series of lemmas in order to to enumerate the preimages under the~$D$ map.  Proposition \ref{prop-D-1} below (and its proof) shows that for any $\al\in\De(\M_n)$, the idempotents from $D^{-1}(\al)$ are obtained from $\al$ by selecting some collection $\th_1,\ldots,\th_t$ of cycle components from $\Ga_\al$, each containing at least one upper outer hook and at least one lower outer hook, and then replacing $2t$ such hooks (one upper and one lower outer hook from each component) by suitable transversals.  The number of idempotents in $D^{-1}(\al)$ corresponding to the collection $\th_1,\ldots,\th_t$ is found by calculating the numbers of upper and lower outer hooks of these components and multiplying all $2t$ of these numbers together.  The total size of~$D^{-1}(\al)$ is then the sum of all such products over all collections of cycle components; algebraically, this sum of products then may be simplified into a single product.
For the statement of the next result, if~$0\leq r\leq n$, we will write $E_r(\M_n)=\set{\alpha\in E(\M_n)}{\rank(\alpha)=r}$.

\begin{prop}\label{prop-D-1}
Let $\alpha\in \De(\M_n)$, and write $q=\rank(\alpha)$ and $k=|\Theta(\alpha)|$.  
\begin{enumerate}[\rm (i)]
\item For any $\beta\in D^{-1}(\alpha)$, $\rank(\beta)=q+2t$ for some $0\leq t\leq k$.
\item For any $0\leq t\leq k$, the set $D^{-1}(\alpha)\cap E_{q+2t}(\M_n)$ has cardinality $\displaystyle{\sum_{\Psi\subseteq\Theta(\alpha)\atop|\Psi|=t} \prod_{\theta\in\Psi} ( u_\theta(\alpha)l_\theta(\alpha))}$.
\vspace{-6mm}
\item We have $|D^{-1}(\alpha)| = \displaystyle{\prod_{\theta\in\Theta(\alpha)}(u_\theta(\alpha)l_\theta(\alpha)+1)}$.
\end{enumerate}
\end{prop}

\begin{proof}
(i).  Suppose $\beta\in D^{-1}(\alpha)$.  Let $t\geq0$ be minimal so that $\alpha=d^t(\be)$.  In the sequence
\[
\beta,d(\beta),\ldots,d^t(\beta)=\alpha,
\]
the rank of each term is $2$ more than the rank of the next term, by definition of the $d$ map.
It follows that $\rank(\beta)=\rank(\alpha)+2t=q+2t$.
We have already noted that $t\geq0$.  By Lemma \ref{lem-d},
\[
\Theta(\beta) \subsetneq \Theta(d(\beta)) \subsetneq\cdots\subsetneq \Theta(d^t(\beta)) = \Theta(\alpha).
\]
Thus, $k=|\Theta(\alpha)|\geq|\Theta(\beta)|+t\geq t$.

\pfitem{ii}  Fix some $0\leq t\leq k$, and write
\begin{align*}
\Sigma=D^{-1}(\alpha)\cap E_{q+2t}(\M_n)
&\AND \sigma=\sum_{\Psi\subseteq\Theta(\alpha)\atop|\Psi|=t} \prod_{\theta\in\Psi} ( u_\theta(\alpha)l_\theta(\alpha)).
\intertext{If $t=0$, then $\Sigma=\{\alpha\}$ and $\sigma=\prod_{\theta\in \emptyset} ( u_\theta(\alpha)l_\theta(\alpha))=1$, as the latter is an empty product.  Now suppose~$t\geq1$.  To complete the proof of (ii), it suffices to find~mutually~inverse~maps}
f:\Sigma\to\bigcup_{\Psi\subseteq\Theta(\alpha)\atop|\Psi|=t} \prod_{\theta\in\Psi} ( U_\theta(\alpha)\times L_\theta(\alpha)) 
&\AND
g:\bigcup_{\Psi\subseteq\Theta(\alpha)\atop|\Psi|=t} \prod_{\theta\in\Psi} ( U_\theta(\alpha)\times L_\theta(\alpha))  \to\Sigma.
\end{align*}
Here, ``$\prod_{\theta\in\Psi}$'' denotes the direct product.

To define $f$, let $\beta\in\Sigma$, and write $\dom(\beta)=\{i_1,\ldots,i_{q+2t}\}$ and $\codom(\beta)=\{l_1,\ldots,l_{q+2t}\}$, where $i_1<\cdots<i_{q+2t}$ and $l_1<\cdots<l_{q+2t}$.  By $t$ applications of Lemma \ref{lem-d}, we see that for each $1\leq h\leq t$, the pair $\big(\{i_{2h-1},i_{2h}\},\{l_{2h-1}',l_{2h}'\}\big)$ belongs to $U_{\th_h}(\al)\times L_{\th_h}(\al)$ for some $\th_h\in\Th(\al)$, and that the components $\th_1,\ldots,\th_t\in\Th(\al)$ are distinct.  So we may define
\[
f(\beta) = \Big( \big( \{i_1,i_2\} , \{l_1',l_2'\} \big) , \ldots , \big(\{ i_{2t-1},i_{2t}\},\{l_{2t-1}',l_{2t}'\}\big) \Big).
\]

To define $g$, let $\Psi=\{\theta_1,\ldots,\theta_t\}\subseteq\Theta(\alpha)$ with $\min(\theta_1)<\cdots<\min(\theta_t)$ and, for each ${1\leq h\leq t}$, let $\{i_{2h-1},i_{2h}\}\in U_{\theta_h}(\alpha)$ and $\{l_{2h-1}',l_{2h}'\}\in L_{\theta_h}(\alpha)$, where $i_{2h-1}<i_{2h}$ and $l_{2h-1}<l_{2h}$.  By Lemma~\ref{lem-upper}, it follows that $i_1<\cdots<i_{2t}$ and $l_1<\cdots<l_{2t}$.  Since these vertices belong to unnested edges of $\alpha$, we may define $\beta\in\M_n$ to be the Motzkin element obtained from $\alpha$ by replacing the non-transversals $\{i_1,i_2\} ,  \ldots , \{ i_{2t-1},i_{2t}\}$ and $\{l_1',l_2'\} , \ldots ,\{l_{2t-1}',l_{2t}'\}$ by the transversals $\{i_1,l_1'\},\ldots,\{i_{2t},l_{2t}'\}$.  Note that $\beta$ is obtained from $\alpha$ by $t$ applications of the process described in Lemma \ref{lem-d-1}, treating the components in the order $\theta_t,\ldots,\theta_1$.  In particular, $\beta\in E(\M_n)$ and $\alpha=D(\beta)$.  By construction, $\rank(\beta)=\rank(\alpha)+2t=q+2t$.  It follows that $\beta\in\Sigma$, so we may then define
\[
g  \Big( \big( \{i_1,i_2\} , \{l_1',l_2'\} \big) , \ldots , \big( \{i_{2t-1},i_{2t}\},\{l_{2t-1}',l_{2t}'\}\big)  \Big) = \beta.
\]
It is easy to check that $f$ and $g$ are mutual inverses.

\pfitem{iii}  We use parts (i) and (ii), and the identity $\prod_{i\in I}(a_i+1)=\sum_{J\sub I}\prod_{j\in J}a_j$, to calculate
\[
|D^{-1}(\alpha)| 
= \sum_{t=0}^k \sum_{\Psi\subseteq\Theta(\alpha)\atop|\Psi|=t} \prod_{\theta\in\Psi} ( u_\theta(\alpha)l_\theta(\alpha))
= \sum_{\Psi\subseteq\Theta(\alpha)} \prod_{\theta\in\Psi} ( u_\theta(\alpha)l_\theta(\alpha))
= \prod_{\theta\in\Theta(\alpha)}(u_\theta(\alpha)l_\theta(\alpha)+1).
\qedhere
\]
\end{proof}

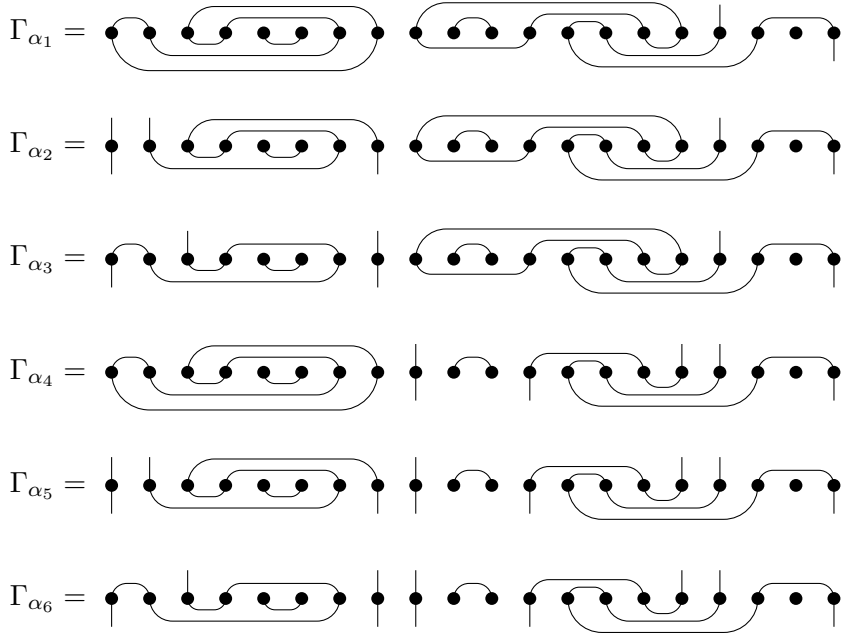
\begin{figure}[h]
\begin{center}
\begin{tikzpicture}[scale=.5]

\begin{scope}[shift={(0,0)}]	
\duarcx34{.3}
\duarcx56{.3}
\duarc9{12}
\duarc{15}{16}
\duarcx27{.6}
\duarcx{14}{17}{.6}
\duarcx{13}{18}{.9}
\darc12
\darc47
\darc{10}{11}
\darc{18}{20}
\darcx{13}{14}{.3}
\darcx{12}{15}{.5}
\darcx9{16}{.8}
\duarcx181
\darcx38{.7}
\stlined{20}
\stlineu{17}
\draw(0.6,0)node[left]{$\Gamma_{\al_1}=$};
\lvs{1,...,20}
\end{scope}

\begin{scope}[shift={(0,-3)}]	
\duarcx34{.3}
\duarcx56{.3}
\duarc9{12}
\duarc{15}{16}
\duarcx27{.6}
\duarcx{14}{17}{.6}
\duarcx{13}{18}{.9}
\darc47
\darc{10}{11}
\darc{18}{20}
\darcx{13}{14}{.3}
\darcx{12}{15}{.5}
\darcx9{16}{.8}
\darcx38{.7}
\stlined{20}
\stlineu{17}
\draw(0.6,0)node[left]{$\Gamma_{\al_2}=$};
\lvs{1,...,20}
\stlineu1
\stlineu2
\stlined1
\stlined8
\end{scope}

\begin{scope}[shift={(0,-6)}]	
\duarcx34{.3}
\duarcx56{.3}
\duarc9{12}
\duarc{15}{16}
\duarcx27{.6}
\duarcx{14}{17}{.6}
\duarcx{13}{18}{.9}
\darc12
\darc47
\darc{10}{11}
\darc{18}{20}
\darcx{13}{14}{.3}
\darcx{12}{15}{.5}
\darcx9{16}{.8}
\stlined{20}
\stlineu{17}
\draw(0.6,0)node[left]{$\Gamma_{\al_3}=$};
\lvs{1,...,20}
\stlineu3
\stlineu8
\stlined1
\stlined8
\end{scope}

\begin{scope}[shift={(0,-9)}]	
\duarcx34{.3}
\duarcx56{.3}
\duarc{15}{16}
\duarcx27{.6}
\duarcx{14}{17}{.6}
\duarcx{13}{18}{.9}
\darc12
\darc47
\darc{10}{11}
\darc{18}{20}
\darcx{13}{14}{.3}
\darcx{12}{15}{.5}
\duarcx181
\darcx38{.7}
\stlined{20}
\stlineu{17}
\draw(0.6,0)node[left]{$\Gamma_{\al_4}=$};
\lvs{1,...,20}
\stlineu9
\stlineu{16}
\stlined9
\stlined{12}
\end{scope}

\begin{scope}[shift={(0,-12)}]	
\duarcx34{.3}
\duarcx56{.3}
\duarc{15}{16}
\duarcx27{.6}
\duarcx{14}{17}{.6}
\duarcx{13}{18}{.9}
\darc47
\darc{10}{11}
\darc{18}{20}
\darcx{13}{14}{.3}
\darcx{12}{15}{.5}
\darcx38{.7}
\stlined{20}
\stlineu{17}
\draw(0.6,0)node[left]{$\Gamma_{\al_5}=$};
\lvs{1,...,20}
\stlineu1
\stlineu2
\stlined1
\stlined8
\stlineu9
\stlineu{16}
\stlined9
\stlined{12}
\end{scope}

\begin{scope}[shift={(0,-15)}]	
\duarcx34{.3}
\duarcx56{.3}
\duarc{15}{16}
\duarcx27{.6}
\duarcx{14}{17}{.6}
\duarcx{13}{18}{.9}
\darc12
\darc47
\darc{10}{11}
\darc{18}{20}
\darcx{13}{14}{.3}
\darcx{12}{15}{.5}
\stlined{20}
\stlineu{17}
\draw(0.6,0)node[left]{$\Gamma_{\al_6}=$};
\lvs{1,...,20}
\stlineu3
\stlineu8
\stlined1
\stlined8
\stlineu9
\stlineu{16}
\stlined9
\stlined{12}
\end{scope}

\end{tikzpicture}
\end{center}
\vspace{-5mm}
\caption{Interface graphs of the six Motzkin idempotents $\al_1,\ldots,\al_6\in E(\M_{20})$ satisfying $D(\al)=\delta$, where $\delta\in\Delta(\M_{20})$ is pictured in Figure \ref{fig:M20cont}.  Note that $d(\al_1)=d(\al_2)=d(\al_3)=d(\al_4)=\al_1=\delta$, while ${d(\al_5)=d(\al_6)=\al_4}$.  Note also that $\al_3$ is the Motzkin element $\gamma$ from Figure \ref{fig:M20}.}
\label{fig:M20contcont}
\end{figure}

To continue the example started above, let $\delta\in\Delta(\M_{20})$ be as in Figure \ref{fig:M20cont}.  Using the values calculated before the statement of Lemma \ref{lem-d-1}, Proposition \ref{prop-D-1}(iii) gives
\[
|D^{-1}(\delta)| = (u_{\theta_1}(\delta)l_{\theta_1}(\delta)+1)(u_{\theta_2}(\delta)l_{\theta_2}(\delta)+1) = (1\cdot2+1)(1\cdot1+1)=6.
\]
The interface graphs of the six elements of $D^{-1}(\delta)$ are depicted in Figure \ref{fig:M20contcont}.

Lemma \ref{lem-EMJ} and Proposition \ref{prop-D-1}(iii) immediately give the following.

\begin{thm}\label{thm-EMJ}
If $S$ is one of $\M_n$ or $\J_n$, then $\displaystyle{|E(S)| = \sum_{\alpha\in \De(S)}\prod_{\theta\in\Theta(\alpha)}(u_\theta(\alpha)l_\theta(\alpha)+1)}$. \qed

\end{thm}

To give the corresponding statement for the Kauffman monoid $\K_n$, for $\alpha\in\M_n$, we write $\Xi(\alpha)$ for the set of \emph{all} cycle components of the interface graph $\Gamma_\alpha$, noting that $\Theta(\alpha)\subseteq\Xi(\alpha)$.  We may identify $\J_n$ with a subset (but not a submonoid) of $\K_n$, by identifying $\alpha\in\J_n$ with $(0,\alpha)\in\K_n$.  By Corollary \ref{cor-kauffman-idempotents}, it follows that $E(\K_n)\subseteq E(\J_n)$.

\begin{thm}\label{thm-EK}
We have $\displaystyle{|E(\K_n)| = \sum_{\alpha\in \De(\J_n)\atop\Xi(\alpha)=\Theta(\alpha)}\prod_{\theta\in\Theta(\alpha)}(u_\theta(\alpha)l_\theta(\alpha))}$.  

\end{thm}

\begin{proof}
First, note that
\[
E(\K_n) = E(\J_n) \cap E(\K_n) = \left( \bigcup_{\alpha\in \De(\J_n)} D^{-1}(\alpha) \right) \cap E(\K_n) = \bigcup_{\alpha\in \De(\J_n)} \big(D^{-1}(\alpha)\cap E(\K_n)\big).
\]
As the sets $D^{-1}(\alpha)$, $\al\in\De(\J_n)$, are pairwise disjoint, it follows that $|E(\K_n)| = \sum_{\alpha\in \De(\J_n)}|D^{-1}(\alpha)\cap E(\K_n)|$.  So it remains to show that
\[
|D^{-1}(\alpha)\cap E(\K_n)| = \begin{cases}
\prod_{\theta\in\Theta(\alpha)}(u_\theta(\alpha)l_\theta(\alpha)) &\text{if $\Xi(\alpha)=\Theta(\alpha)$}\\
0 &\text{otherwise.}
\end{cases}
\]
With this in mind, let $\alpha\in \De(\J_n)$, and write $q=\rank(\alpha)$.  If there is a cycle component ${\theta\in\Xi(\alpha)\setminus\Theta(\alpha)}$, then $\theta$ is a cycle component of any $\beta\in D^{-1}(\alpha)$, by Lemma \ref{lem-d-1}, and it then follows from Corollary \ref{cor-kauffman-idempotents} that ${D^{-1}(\alpha)\cap E(\K_n)=\emptyset}$.  Next, suppose $\Xi(\alpha)=\Theta(\alpha)$, and put $k=|\Theta(\alpha)|$.  By Proposition \ref{prop-D-1}(i), an element of $D^{-1}(\alpha)$ has rank $q+2t$ for some $0\leq t\leq k$.  The interface graph of such an element contains $k-t$ cycle components, so (again using Corollary \ref{cor-kauffman-idempotents}) we only obtain an element of $E(\K_n)$ in the case $t=k$, in which case \emph{all} elements of $D^{-1}(\alpha)\cap E_{q+2t}(\M_n)$ belong to $E(\K_n)$.  Proposition \ref{prop-D-1}(ii) then gives the stated value of $|D^{-1}(\alpha)\cap E(\K_n)|=|D^{-1}(\alpha)\cap E_{q+2k}(\M_n)|$, since the only term in the sum in Proposition \ref{prop-D-1}(ii) when $t=k$ is the $\Psi=\Theta(\alpha)$ term.
\end{proof}

\begin{rem}\label{rem:TL2}
As in Remark \ref{rem:TL1}, Theorem \ref{thm-EK} also gives the number of idempotent basis elements of the Temperley-Lieb and Motzkin algebras in the generic case.  As in \cite[Section 6]{DEEFHHL1}, the formula in Theorem \ref{thm-EK} could be adapted to treat the case in which the twisting parameter $\xi$ is an $M$th root of unity with $M\leq n$, but we omit the details.
\end{rem}

We may also give formulae for the number of idempotents of $\M_n,\J_n,\K_n$ of fixed rank.  Recall that for $0\leq r\leq n$, we write $E_r(\M_n)=\set{\alpha\in E(\M_n)}{\rank(\alpha)=r}$.  If $S$ is one of $\J_n$ or~$\K_n$, we will write $E_r(S)=S\cap E_r(\M_n)$.  Recall that we are identifying $E(\K_n)$ with a subset of $E(\J_n)$.  Note that $E_r(\J_n)=E_r(\K_n)=\emptyset$ if $r\not\equiv n\pmod2$.  

\begin{thm}\label{thm-ErMJK}
Let $0\leq r\leq n$, and write $r=q+2t$ where $q\in\{0,1\}$.
\begin{enumerate}[\rm (i)]
\item If $S$ is one of $\M_n$ or $\J_n$, then $|E_r(S)| = \displaystyle{\sum_{\alpha\in \De(S)\atop\rank(\alpha)=q}\sum_{\Psi\subseteq\Theta(\alpha)\atop|\Psi|=t} \prod_{\theta\in\Psi} ( u_\theta(\alpha)l_\theta(\alpha))}$.

\item We have $|E_r(\K_n)| = \displaystyle{\sum_\alpha \prod_{\theta\in\Theta(\alpha)}(u_\theta(\alpha)l_\theta(\alpha))}$,
where the sum is over all $\alpha\in \De(\J_n)$ with $\rank(\alpha)=q$, $\Xi(\alpha)=\Theta(\alpha)$ and $|\Theta(\alpha)|=t$.
\end{enumerate}
\end{thm}

\begin{proof}
Part (i) follows quickly from Proposition \ref{prop-D-1} (and its proof), and part (ii) from the proof of Theorem~\ref{thm-EK}.
\end{proof}

\section{The algorithms}\label{sect-algorithms}

In this section we describe algorithms for enumerating the idempotents in the Jones, Kauffman and Motzkin monoids, based on the theoretical results obtained in Section \ref{sect-theory}.  These are presented in Algorithms \ref{algorithm-jones}--\ref{algorithm-motzkin}, below.

In the algorithms described in this section, it is necessary to enumerate the
interface graphs of the elements of the Jones monoid $\J_{n}$ of minimal rank, and of the
elements of the Motzkin monoid $\M_{n}$ of ranks 0 and 1. Roughly speaking, in accordance with
Theorems~\ref{thm-EMJ} and \ref{thm-EK}, the algorithms then involve finding connected components of
these interface graphs and counting the number of upper and lower outer hooks involved in every
cycle component. This could be achieved using standard graph theoretic
algorithms, and there would be essentially nothing further to describe. In
general, however, determining the connected components of a graph with $v$ vertices and
$e$ edges has complexity $O(v+e)$; see \cite{Tarjan1972aa,Gabow2000aa,Sedgewick1983aa}. Hence, in the case of the even degree
Jones monoid $\J_{2n}$, for example, the complexity of this approach would be
$O(4nC_n^2)$, where $C_n=\frac1{n+1}\binom{2n}n
$ is the $n$th Catalan number.  This is because the interface graph of any element of $\De(\J_{2n})$ has $2n$ vertices and $2n$ edges, and since $\De(\J_{2n})$ has size $C_n^2$ \cite[Proposition 2.7(iii)]{DEG2017}.
However, the interface graphs under consideration have several special
properties that we can exploit to substantially reduce the run time of our algorithms.
In order to explain these properties, and hence the nature of the algorithms, we must first discuss a number of concepts relevant to the enumeration of interface graphs.

\subsection{Background on Dyck and Motzkin words}

A \textit{Dyck word} is a balanced string of left and right brackets.
\textit{Balanced} in this context means that the numbers of left and right
brackets are equal, and that at any point, reading from left to right, the number of
left brackets is at least the number of right brackets. A Dyck word necessarily has even length.  For example, $u=((())())$ is a Dyck word of length $8$.  
A recent algorithm of Neri \cite{Neri2016aa,Neri2016ab} allows for fast generation of Dyck words.
%
We will write $D_{2n}$ for the set of all Dyck words of length $2n$.  So $|D_{2n}|=C_n$ is the $n$th Catalan number.

A \textit{Motzkin word} is a string of left and right brackets and dots,
such that the subword consisting only of the brackets is a Dyck word.  A Motzkin word can have any length.  For example, $v=()(\cdot()){\cdot}(\cdot\cdot)(()())$ is a Motzkin word of length $18$.  A Motzkin word of length $n$ can be thought of as a pair consisting of a Dyck word of length $m$, for some even $m\leq n$, and a subset of $\bn=\{1,\ldots,n\}$ of size $n-m$; the subset specifies the positions of the dots in the Motzkin word.  For example, the above Motzkin word $v$ corresponds to the Dyck word $()(())()(()())\in D_{14}$ and the subset $\{4, 8, 10, 11\}$.  
We will write $M_n$ for the set of all Motzkin words of length $n$; there should be no confusion with the Motzkin monoid itself, which is denoted by $\M_n$.  So $|M_n|=\mu(n,0)$, where the numbers $\mu(n,r)$ were defined in \eqref{equation-motzkin1} and \eqref{equation-motzkin2}.
It is relatively straightforward to produce the sets $M_n$ for the values of $n$ we are concerned with here,
namely for $n\leq20$. For instance, there are 50\ 852\ 019 Motzkin
words of length~20, and these can be produced and stored in a convenient format for
use in Algorithm~\ref{algorithm-motzkin} in about~35 seconds, and using about 2
GB of memory.  In particular, creating and storing the Motzkin words of a given length represents a tiny fraction of the time taken by Algorithm \ref{algorithm-motzkin}.  We will not describe the process for producing the Motzkin words in more detail here.

Denote by $\S_n$ the \emph{symmetric group} of degree $n$, which consists of all permutations of~$\bn$.  There is a natural injective map
\[
\bp: M_n \to \S_n
\]
taking a Motzkin word $w\in M_n$ to the permutation $\bp(w)\in\S_n$ defined, for $i\in\bn$, by
\[
\bp(w)(i) = \begin{cases}
i &\text{if $w$ has a dot in position $i$}\\
j &\text{if $w$ has a bracket in position $i$ that is matched by a bracket in position $j$.}
\end{cases}
\]
So, if we suppose $w\in M_n$ has left brackets at positions $i_1<\cdots<i_m$, and that these are matched by right brackets at positions $j_1,\ldots,j_m$, respectively, then $\bp(w)$ can be written as a product of commuting transpositions as $\bp(w)=(i_1\ j_1)\cdots(i_m\ j_m)$.  For example, with $u\in D_8$ and $v\in M_{18}$ as defined above,
\[
\bp(u)=(1\ 8)(2\ 5)(3\ 4)(6\ 7)\in\S_8 \ANd \bp(v)=(1\ 2)(3\ 7)(5\ 6)(9\ 12)(13\ 18)(14\ 15)(16\ 17)\in\S_{18}.
\]
In general, for a Motzkin word $w\in M_n$, $\bp(w)$ has no fixed points if and only if~$w\in D_n$; in this case, $n$ must be even.  Note also that if $w\in M_n$ has at least one bracket, then $\bp(w)$ is an involution (a permutation of order $2$); if $w$ consists only of dots, then $\bp(w)$ is the identity element of $\S_n$.

Now consider a Motzkin element $\al\in\M_n$ with $\rank(\al)=0$.  The upper blocks of $\al$ induce a Motzkin word $w_1\in M_n$ in a natural way; for each upper hook $\{i,j\}$ of $\al$ with $i<j$, $w_1$ has a left bracket at position~$i$ and a right bracket at position $j$, while the upper singletons of~$\al$ correspond to the dots of $w_1$.  Similarly, the lower blocks of $\al$ induce a second Motzkin word $w_2\in M_n$.  We will write $\bm(\al)=(w_1,w_2)\in M_n\times M_n$ for the pair consisting of these two words.  Conversely, given any pair $(w_1,w_2)\in M_n\times M_n$, it is easy to see that there is a Motzkin element $\al\in\M_n$ of rank $0$ with $\bm(\al)=(w_1,w_2)$.

To describe the Motzkin elements of rank $1$ in terms of Motzkin words, we first define $M_{n+1}'$ to be the subset of $M_{n+1}$ consisting of all Motzkin words of length $n+1$ whose last symbol is a right bracket.  Equivalently, a Motzkin word $w\in M_{n+1}$ belongs to $M_{n+1}'$ if $\bp(w)(n+1)\not=n+1$.  Now consider a Motzkin element $\al\in\M_n$ with $\rank(\al)=1$, and let the unique transversal of $\al$ be $\{i,j\}$.  We define the Motzkin word $w_1\in M_{n+1}'$ to have a left bracket at position $i$, a right bracket at position $n+1$, and where the remaining symbols of $w_1$ are determined by the upper blocks of $\al$ in the same way as in the previous paragraph.  We define $w_2$ analogously, in terms of $j$ and the lower blocks of $\al$.  Again, we will write ${\bm(\al)=(w_1,w_2)\in M_{n+1}'\times M_{n+1}'}$ for the pair consisting of these two words.  Again, given any pair ${(w_1,w_2)\in M_{n+1}'\times M_{n+1}'}$, there is a Motzkin element $\al\in\M_n$ of rank $1$ with $\bm(\al)=(w_1,w_2)$.

The previous two paragraphs describe a bijection
\[
\bm : \set{\al\in\M_n}{\rank(\al)\leq1} \to (M_n\times M_n) \cup (M_{n+1}'\times M_{n+1}').
\]
%
%
We denote by $\bd$ the restriction of $\bm$ to the Jones elements $\set{\al\in\J_n}{\rank(\al)\leq1}$ of rank at most~$1$.  The restriction $\bd$ is a bijection onto its image, which is $D_{2\lceil\frac n2\rceil}\times D_{2\lceil\frac n2\rceil}$: that is, either $D_n\times D_n$ or $D_{n+1}\times D_{n+1}$, according to whether $n$ is even or odd, respectively.  In particular, Jones elements of minimum rank correspond to certain pairs of Dyck words of an appropriate length.


We noted above that $|D_{2n}|=C_n$ and that $|M_n|=\mu(n,0)$ for any $n$.  It is also known \cite[Proposition~2.8]{DEG2017} that ${|M_{n+1}'|=\mu(n,1)}$.
%
%
%
In the algorithms presented in this section we will fix (arbitrary) orderings on the sets $D_{2n}$, $M_n$ and $M_{n+1}'$, and will denote the elements of these sets as
\[
D_{2n} = \set{u_i}{1\leq i \leq C_n},\quad
M_n = \set{w_i}{1\leq i\leq \mu(n, 0)\}},\quad
M_{n+1}' = \set{w_i'}{1 \leq i\leq \mu(n, 1)}.
\]
If $w\in M_n$ is a Motzkin word, then we say that a left bracket of $w$ is an \textit{outer
bracket} if this bracket is not enclosed by any other brackets.  We define $\bo(w)$ to be the subset of $\bn$ for which $i\in\bo(w)$ if and only if~$w$ has an outer bracket at position $i$.
For example, for $u\in D_8$ and $v\in M_{18}$ defined above, $\bo(u)=\{1\}$ and~$\bo(v)=\{1,3,9,13\}$.

Finally, recall that in Subsection \ref{sect-mapping} we defined and studied a map
\[
D:E(\M_n)\to\set{\al\in E(\M_n)}{\rank(\al)\leq1}.
\]
In this section, for convenience, we will write $\widehat\al=D(\al)$ for any $\al\in E(\M_n)$.


\subsection{The algorithm for Jones idempotents}\label{subsect-Jones}

Algorithm~\ref{algorithm-jones} contains pseudocode for counting the number of
idempotents in the Jones monoid~$\J_{n}$.  A  C++ implementation of this algorithm can be found at
\cite{Mitchell2016aa}.  
Roughly speaking, the algorithm begins by enumerating the elements of $\De(\J_n)$ in terms of pairs $(u_i,u_j)$ of Dyck words of length $n$ or $n+1$, as appropriate.  It then proceeds to count the outer hooks in each connected component of the interface graph of the Jones element $\bd^{-1}(u_i,u_j)\in\De(\J_n)$; this then yields 
the number of idempotents $\al\in E(\J_n)$ with~${\bd(\widehat\al)=(u_i,u_j)}$, according to Proposition \ref{prop-D-1}(iii).  The algorithm then concludes by summing these values.

\begin{algorithm}[h]
  \caption{Count the number of idempotents in the Jones monoid $\J_{n}$}
  \label{algorithm-jones}
  \begin{algorithmic}[1]
    \State $N := 0$ \Comment{Number of idempotents}
    \For{$i \in \{1, \ldots, C_{\lceil n/2\rceil}\}$}
    \Comment{Loop over Dyck words of length $n$ or $n + 1$}
    \State $N \gets N + 2 ^ {|\bo(u_i)\setminus \{\bp(u_i)(n +
    1)\}|}$ 
    \label{algorithm-jones-check-0}
    \Comment{$\al\in E(\J_{n})$
    such that $\bd(\widehat{\al}) = (u_i, u_i)$}
    \For{$j \in \{i + 1, \ldots, C_{\lceil n/2\rceil}\}$} 
    \Comment{Loop over Dyck words}
    \State $M := 1$ \Comment{Number of idempotents $\al$ with
    $\bd(\widehat{\al}) = (u_i, u_j)$}
    \State $m := 0$ \Comment{Largest value seen in any cycle of
    $\Gamma_{\bd^{-1}(u_i, u_j)}$}
    \While{$m < \max\{\bo(u_j)\setminus \{\bp(u_j)(n +
    1)\}\}$} 
    \label{algorithm-jones-check-00}
    \Comment{Loop over cycles of
    $\Gamma_{\bd^{-1}(u_i, u_j)}$}
    \State $k, l \gets \min\set{x\in \bo(u_j)}{x \geq m}$
    \Comment{Start of the next cycle}
    \State $I, J := 0$ 
    \Comment{Count the number of outer hooks in this cycle}
    \Repeat\Comment{Loop within the current cycle}
    \If{$l \in \bo(u_i)$}\label{algorithm-jones-check-1} \Comment{Found
    an outer bracket of $u_i$ in current cycle}
    \State $I \gets I + 1$
    \EndIf
    \If{$l \in \bo(u_j)$}\label{algorithm-jones-check-2} \Comment{Found
    an outer bracket of $u_j$ in current cycle}
    \State $J \gets J + 1$
    \State\label{testing} $m \gets \max\{m,\ \bp(u_j)(l)\}$
    \EndIf
    \State $l \gets \bp(u_i)\bp(u_j)(l)$ \Comment{Go to the next
    position in the current cycle}
    \Until{{$l = k$}} \Comment{Returned to the start of the cycle}
    \State $M \gets M(IJ + 1)$ \Comment{Multiply by number of outer brackets in
    current cycle}
    \EndWhile
    \State $N \gets N + 2M$ \Comment{Add number of idempotents $\al$ with
    $\bd(\widehat{\al}) = (u_i, u_j)$}
    \EndFor
    \EndFor
    \State \Return $N$
  \end{algorithmic}
\end{algorithm}

Before moving on to the other algorithms, we first comment on a number of features of Algorithm \ref{algorithm-jones}, including some simple optimisations that have been included.  

First, if $\al\in \De(\J_n)$ is such that $\bd(\al)=(u_i,u_i)$ for some $i$, then every component of the interface graph~$\Ga_\al$ is a cycle of length $2$ or an active path of length $0$ (the latter only occurs when $n$ is odd, in which case there is a unique such path).  As such, Proposition \ref{prop-D-1}(iii) tells us that $D^{-1}(\al)
$ has size $2^k$ if $n$ is even, or $2^{k-1}$ if~$n$ is odd, respectively, where $k$ is the number of outer (left) brackets of $u_i$.  The reason for subtracting~$1$ from $k$ in the case~$n$ is odd is that the last outer bracket of $u_i\in D_{n+1}$ corresponds to the path component of $\Ga_\al$.  See Line \ref{algorithm-jones-check-0} of Algorithm \ref{algorithm-jones}.  
Lines \ref{algorithm-jones-check-0} and \ref{algorithm-jones-check-00} refer to $\bp(u_i)(n + 1)$, which is only defined when $n$ is odd, and which can be ignored when $n$ is even.
In the implementation in \cite{Mitchell2016aa}, Algorithm \ref{algorithm-jones} is split into two parts covering the even and odd cases separately.


If $u_i$ and $u_j$ are distinct Dyck words, then there are the same number of idempotents $\al\in E(\J_n)$ with~${\bd(\widehat\al)=(u_i,u_j)}$ as there are with $\bd(\widehat\al)=(u_j,u_i)$.  This corresponds to the anti-involution ${}^*:\J_n\to\J_n$, and the fact that a Jones element $\al$ is an idempotent if and only if $\al^*$ is, since $[\bd^{-1}(u_i,u_j)]^*=\bd^{-1}(u_j,u_i)$.   

A further optimization along these lines is available in the case that $n$ is even.  Namely, for any~$n$, whether even or odd, there is an involution ${}^\dagger: \J_n \to \J_n$, where $\al^\dagger$ is the result of reflecting~$\al$ in a vertical axis midway between points $1$ and $n$.  The involution ${}^\dagger$ was studied along with the anti-involution ${}^*$ in \cite{ADV2012_2}.  It is again clear that $\al\in\J_n$ is an idempotent if and only if $\al^\dagger$ is.  For an even value of $n$, and for $\al\in\J_n$ of rank~$0$, if $\bd(\al)=(u_i,u_j)$, then $\bd(\al^\dagger)=(\operatorname{rev}(u_i),\operatorname{rev}(u_j))$, where $\operatorname{rev}(w)$ is the result of writing $w$ in reverse and interchanging left and right brackets.  This means that in the case that $\al$ and $\al^\dagger$ are not equal, we only need to calculate the size of one of $D^{-1}(\al)$ or $D^{-1}(\al^\dagger)$.  The implementation of this optimisation is rather technical, and only applies in the even case (since the active paths in the interface graphs of rank $1$ Jones elements interfere with the ${}^\dagger$ map in the case of odd $n$), so we have not included it in the pseudocode for Algorithm~\ref{algorithm-jones}.  This optimisation is included in the implementation~\cite{Mitchell2016aa}.

As a further note, it is not necessary to check if $l\in \bo(u_i)$ and
$l\in \bo(u_j)$ in Lines~\ref{algorithm-jones-check-1}
and~\ref{algorithm-jones-check-2} of Algorithm~\ref{algorithm-jones}, since it can be shown that if $l\in \bo(u_i)$, then $l\not\in
\bo(u_j)$, and vice versa, unless $l$ is the minimum vertex in its component. In fact, apart from the above-mentioned exception, it is only possible to have $l\in
\bo(u_j)$ before the first time that $l\in \bo(u_i)$, and this
could be separated into another loop to reduce the number of branches in the
innermost loops.  However, for the sake of brevity we do not include this
optimization in the pseudocode in Algorithm~\ref{algorithm-jones}, although it is
included in the implementation~\cite{Mitchell2016aa}.

Finally, we note that Algorithm~\ref{algorithm-jones} is \emph{embarrassingly parallel}, in
the sense that the number of idempotents $\al\in E(\J_{n})$ such that $\bd(\widehat\al)
= (u_i, u_j)\in D_{2\lceil n/2\rceil }\times D_{2\lceil n/2\rceil}$ can be
enumerated independently for different values of $i$ and $j$. 
%


\subsection{The algorithm for Kauffman idempotents}

The key difference between Algorithm~\ref{algorithm-jones} for the Jones monoid $\J_n$
and Algorithm~\ref{algorithm-kauffman} (below) for the Kauffman monoid $\K_n$ is that if, for some $\al\in \De(\J_n)$,
there is a cycle in the interface graph $\Gamma_\al$
containing no upper outer hooks or no lower outer hooks, then $D^{-1}(\al)\cap E(\K_n)$ is empty; see the proof of Theorem \ref{thm-EK}.
Hence, in
Algorithm~\ref{algorithm-kauffman}, \emph{every} cycle of $\Gamma_\al$ must be considered and not only those starting at an outer hook as was
the case in Algorithm~\ref{algorithm-jones}.  Again, a number of optimisations are included in Algorithm~\ref{algorithm-kauffman}, but we will not comment explicitly on these, as they are virtually identical to those in Algorithm~\ref{algorithm-jones}. 

Note that for any $1\leq i\leq C_{\lceil n/2\rceil}$, there are no idempotents $\al\in E(\K_n)$ with $\bd(\widehat\al)=(u_i,u_i)$ unless $u_i=()()()\ldots()$, in which case the identity element is the only such idempotent.  This is why we start with $N=1$ in Line~\ref{line-1-kauffman} of Algorithm \ref{algorithm-kauffman}, and why we only consider pairs $(u_i,u_j)$ with $i\not=j$; see Lines \ref{line-2-kauffman} and \ref{line-3-kauffman}.


\begin{algorithm}[H]
  \caption{Count the number of idempotents in the Kauffman monoid $\K_{n}$}
  \label{algorithm-kauffman}
  \begin{algorithmic}[1]
    \State $N := 1$ \Comment{Number of idempotents}
    \label{line-1-kauffman}
    \For{$i \in \{1, \ldots, C_{\lceil n/2\rceil}\}$}
    \label{line-2-kauffman}
    \Comment{Loop over Dyck words of length $n$ or $n + 1$}
    \For{$j \in \{i + 1, \ldots, C_{\lceil n/2\rceil}\}$} 
    \label{line-3-kauffman}
    \Comment{Loop over Dyck words}
    \State $M := 1$ \Comment{Number of idempotents $\al\in E(\K_n)$ with
    $\bd(\widehat{\al}) = (u_i, u_j)$}
    \State $b := 1$ \Comment{The current vertex}
    \State $B := \varnothing$ 
    \Comment{The vertices of $\Gamma_{\bd^{-1}(u_i, u_j)}$ seen already}
    \While{$M \not= 0$ and $b < n + 1$} 
    \Comment{Loop over cycles of $\Gamma_{\bd^{-1}(u_i, u_j)}$}
    \State $I, J := 0$ 
    \Comment{Count the number of outer hooks in this cycle}
    \Repeat
    \If{$b \in \bo(u_i)$} \Comment{Found an outer bracket of $u_i$ in
    current cycle}
    \State $I \gets I + 1$
    \EndIf
    \If{$b \in \bo(u_j)$} \Comment{Found an outer bracket of $u_j$ in
    current cycle}
    \State $J \gets J + 1$
    \EndIf
    \State $B\gets B\cup \{b\}$
    \State $b \gets \bp(u_i)\bp(u_j)(b)$ \Comment{Go to the next
    position in the current cycle}
    \Until{$b \in B$} 
    \State $M \gets MIJ$ \Comment{Multiply by number of outer brackets in
    current cycle}
    \While{$M \not = 0$ and $b \in B$} 
    \Comment{Find the next vertex $b$ not seen already}
    \State $b\gets b + 1$ 
    \EndWhile
    \EndWhile
    \State $N \gets N + 2M$ \Comment{Add number of idempotents $\al$ with
    $\bd(\widehat{\al}) = (u_i, u_j)$}
    \EndFor
    \EndFor
    \State \Return $N$
  \end{algorithmic}
\end{algorithm}

\subsection{The algorithm for Motzkin idempotents}

Algorithm \ref{algorithm-motzkin} contains pseudocode for calculating the number of idempotents in the Motzkin monoid $\M_n$.  The basic idea of Algorithm \ref{algorithm-motzkin} is similar to Algorithm \ref{algorithm-jones}, except that separate parts are required to count idempotents of even rank (Lines \ref{motz-2}--\ref{motz-20}) and odd rank (Lines \ref{motz-21}--\ref{motz-41}).  Similar optimisations to Algorithms \ref{algorithm-jones} and \ref{algorithm-kauffman} have been included.

\begin{algorithm}[H]
  \caption{Count the number of idempotents in the Motzkin monoid $\M_{n}$}
  \label{algorithm-motzkin}
  \begin{algorithmic}[1]
    \State $N := 0$ \Comment{Number of idempotents}
    \For{$i \in \{1, \ldots, \mu(n, 0)\}$}
    \label{motz-2}
    \Comment{Loop over Motzkin words of length $n$}
    \State $N \gets N + 2 ^ {|\bo(w_i)|}$ \Comment{$\al\in E(\M_{n})$
    with $\bm(\widehat{\al}) = (w_i, w_i)$}
    \For{$j \in \{i + 1, \ldots, \mu(n, 0)\}$} 
    \Comment{Loop over Motzkin words}
    \State $M := 1$ 
    \Comment{Number of idempotents $\al$ with $\bm(\widehat{\al}) = (w_i,
    w_j)$}
    \State $m := 0$ \Comment{Largest value seen in any path of
    $\Gamma_{\bm^{-1}(w_i, w_j)}$}
    \While{$m < \max\{\bo(w_j)\}$} 
    \Comment{Loop over paths of $\Gamma_{\bm^{-1}(w_i, w_j)}$}
    \State $k, l\gets \min\set{x\in \bo(w_j)}{x \geq m}$
    \Comment{Start of the next path}
    \State $I, J := 0$ 
    \Comment{Count the number of outer hooks in this cycle}
    \Repeat \Comment{Follow the current path}
    \If{$l \in \bo(w_i)$} \Comment{Found an outer bracket of $w_i$ in
    current path}
    \State $I \gets I + 1$
    \EndIf
    \If{$l \in \bo(w_j)$} \Comment{Found an outer bracket of $w_j$ in
    current path}
    \State $J \gets J + 1$
    \State $m \gets \max\{m,\ \bp(w_j)(l)\}$
    \EndIf
    \State $l \gets \bp(w_i)\bp(w_j)(l)$ \Comment{Go to the next
    position in the current path}
    \Until{$l = k$ or $\bp(w_i)(l) = l$ or
    $\bp(w_j)\bp(w_i)(l) = \bp(w_i)(l)$}
    \If{$l = k$} \Comment{The current path is a cycle} 
      \State $M \gets M(IJ + 1)$ \Comment{Multiply by number of outer brackets in
      current cycle}
    \EndIf
    \EndWhile
    \State $N \gets N + 2M$ \Comment{Add number of idempotents $\al$ with
    $\bm(\widehat{\al}) = (w_i, w_j)$}
    \label{motz-20}
    \EndFor
    \EndFor
    
    \For{$i \in \{1, \ldots, \mu(n, 1)\}$}
    \label{motz-21}
    \Comment{Loop over Motzkin words of length $n + 1$ in $M_{n + 1}'$}
    \State $N \gets N + 2 ^ {|\bo(w_i)|}$ \Comment{$\al\in E(\M_{n})$
    with $\bm(\widehat{\al}) = (w_i', w_i')$}
    \For{$j \in \{i + 1, \ldots, \mu(n, 1)\}$} 
    \Comment{Loop over elements of $M_{n + 1}'$}
    \State $M := 1$ 
    \Comment{Number of idempotents $\al$ with $\bm(\widehat{\al}) = (w_i',
    w_j')$}
    \State $m := n + 1$ \Comment{Smallest value seen in any path of
    $\Gamma_{\bm^{-1}(w_i', w_j')}$}
    \While{$m > \min\{\bo(w_j')\}$ and $M \not= 0$} 
    \Comment{Loop over paths of $\Gamma_{\bm^{-1}(w_i', w_j')}$}
    \State $k, l\gets \max\set{\bp(w_j')(x)}{x \in \bo(w_j'),\ x
    \leq m}$
    \Comment{Start of the next path}
    \State $I, J := 0$ 
    \Comment{Count the number of outer hooks in this cycle}
    \Repeat \Comment{Follow the current path}
    \If{$\bp(w_j')(l) \in \bo(w_i')$} 
    \Comment{Found an outer bracket of $w_i'$ in current path}
    \State $I \gets I + 1$
    \EndIf
    \If{$\bp(w_j')(l) \in \bo(w_j')$} 
    \Comment{Found an outer bracket of $w_j'$ in current path}
    \State $J \gets J + 1$
    \State $m \gets \min\{m,\ l\}$
    \EndIf
    \State $l \gets \bp(w_i')\bp(w_j')(l)$ \Comment{Go to the next
    position in the current path}
    \Until{$l = k$ or $\bp(w_i')(l) = l$ or
    $\bp(w_j')\bp(w_i')(l) = \bp(w_i')(l)$}
    \If{$l = k$} \Comment{The current path is a cycle} 
      \State $M \gets M(IJ + 1)$ \Comment{Multiply by number of outer brackets in
      current cycle}
    \ElsIf{$k = n + 1$}\Comment{The current path is not a cycle}  
    \State $M \gets 0$
    \Comment{There are no $\al$ with $\bm(\widehat{\al}) = (w_i', w_j')$ to
    count}
    \EndIf
    \EndWhile
    \State $N \gets N + 2M$ \Comment{Add number of idempotents $\al$ with
    $\bm(\widehat{\al}) = (w_i', w_j')$}
    \label{motz-41}
    \EndFor
    \EndFor
    \State \Return $N$
  \end{algorithmic}
\end{algorithm}

\section{Values and benchmarking}\label{sect-values}

In this section, we give some calculated values of the various number sequences we have considered.  Tables~\ref{tab-EJn-EKn} and~\ref{tab-EMn} give the number of idempotents in the Jones, Kauffman and Motzkin monoids, as well as comparative running times for the various algorithms described in the paper.
For each algorithm, these values were computed using GAP or~\cite{Mitchell2016aa}, as appropriate, on an IBM power8 (8247-22L), with 24 cores at 3,026 GHz (giving 192 threads) running powerKVM.
At the time of writing, these represent the largest known values of $|E(\J_n)|$, $|E(\K_n)|$ and $|E(\M_n)|$; cf.~Sequences A225798, A281438 and A256672 on \cite{OEIS}.
Note that values of $|E(\J_{2n})|$ can be computed faster than $|E(\J_{2n-1})|$ because of the $\dagger$ map discussed in Subsection \ref{subsect-Jones}.

\begin{table}[ht]%
  \begin{center}
  {\small
    \begin{tabular}{|r|r|r|r|r|}
      \hline
      $n$ &    
      \multicolumn{1}{c|}{$|E(\J_n)|$} &
      \multicolumn{1}{c|}{Alg.~\ref{algorithm-basic}} &
      \multicolumn{1}{c|}{Alg.~\ref{algorithm-regular-*}} &
      \multicolumn{1}{c|}{Alg.~\ref{algorithm-jones}} \\
      \hline 
      0   & 1 &&& \\ 
      1  & 1 &&& \\ 
      2  & 2 &&& \\
      3  & 5 &&& \\ 
      4  & 12 &&& \\
      5  & 36 &&& \\ \hline
      6  & 96 &&& \\ 
      7  & 311 &&& \\ 
      8  & 886 &&& \\ 
      9  & 3000 &&& \\ 
      10 & 8944 &&& \\ \hline
      11 & 31\ 192  &&&\\ 
      12 & 96\ 138 &&&  \\ 
      13 & 342\ 562 &2&&\\ 
      14 & 1083\ 028 &8&& \\ 
      15 & 3923\ 351 &32&&\\ \hline
      16 & 12\ 656\ 024&5901&1& \\ 
      17 & 46\ 455\ 770 &-&4&\\ 
      18 & 152\ 325\ 850 &-&16&\\ 
      19 & 565\ 212\ 506 &-&51&\\ 
      20 & 1878\ 551\ 444&-&214& \\ \hline
      21 & 7033\ 866\ 580 &- &689 & 2 \\
      22 & 23\ 645\ 970\ 022 & - &- & 2 \\ 
      23 & 89\ 222\ 991\ 344 & - &- & 29  \\ 
      24 & 302\ 879\ 546\ 290 & - &- & 23  \\ 
      25 & 1150\ 480\ 017\ 950 & - &- & 522  \\ \hline
      26 & 3938\ 480\ 377\ 496 & - & - & 500  \\
      27 & 15\ 047\ 312\ 553\ 918 & - & - & 7260  \\
      28 & 51\ 892\ 071\ 842\ 570 & - & - & 5520  \\
      29 & 199\ 274\ 492\ 098\ 480 & - & - & 101\ 160 \\
      30 & 691\ 680\ 497\ 233\ 180 & - & - & 77\ 100 \\\hline 
    \end{tabular}
\qquad\qquad\qquad
    \begin{tabular}{|r|r|r|}
    \hline
    $n$ &    
    \multicolumn{1}{c|}{$|E(\K_n)|$} &
    \multicolumn{1}{c|}{Alg.~\ref{algorithm-kauffman}} \\ \hline
    0  & 1                  &    \\ 
    1  & 1                  &    \\ 
    2  & 1                  &    \\ 
    3  & 3                  &    \\ 
    4  & 5                  &    \\ 
    5  & 15                 &    \\ \hline
    6  & 31                 &    \\ 
    7  & 93                 &    \\ 
    8  & 215                &    \\ 
    9  & 653                &    \\ 
    10 & 1619             &    \\ \hline
    11 & 4979             &    \\ 
    12 & 12\ 949            &    \\ 
    13 & 40\ 293            &    \\ 
    14 & 108\ 517           &    \\ 
    15 & 341\ 241           &    \\ \hline
    16 & 943\ 937           &    \\ 
    17 & 2996\ 127        &    \\ 
    18 & 8465\ 319        &    \\ 
    19 & 27\ 092\ 419       &    \\ 
    20 & 77\ 878\ 271       &    \\ \hline
    21 & 251\ 073\ 791      & 5     \\ 
    22 & 732\ 129\ 719      & 5     \\ 
    23 & 2375\ 764\ 351   & 60    \\ 
    24 & 7012\ 025\ 277   & 67    \\ 
    25 & 22\ 886\ 955\ 207  & 787   \\ \hline
    26 & 68\ 254\ 122\ 669  & 912   \\ 
    27 & 223\ 946\ 197\ 065 & 10\ 740 \\ 
    28 & 673\ 885\ 100\ 857 & 12\ 300 \\ 
    29 & 2221\ 505\ 541\ 773                  & 147\ 300     \\ 
    30 & 6737\ 598\ 265\ 009                  & 165\ 720     \\ \hline
    \end{tabular}
}
  \end{center}
\vspace{-5mm}
  \caption{
Left: the number of idempotents in the Jones monoid $\J_n$, and the time in seconds to calculate these numbers using Algorithms~\ref{algorithm-basic},~\ref{algorithm-regular-*} and~\ref{algorithm-jones}.
Right: the number of idempotents in the Kauffman monoid $\K_n$, and the time in seconds to calculate these numbers using Algorithm~\ref{algorithm-kauffman}; note that Algorithms \ref{algorithm-basic} and \ref{algorithm-regular-*} do not apply to $\K_n$, as it is neither finite nor a regular $*$-semigroup; see Section \ref{sect-existing}.
}
  \label{tab-EJn-EKn}
\end{table}

  
\begin{table}[ht]%
  \begin{center}
    \begin{tabular}{|r|r|r|r|r|}
    \hline
          $n$ &    \multicolumn{1}{c|}{$|E(\M_n)|$} &
          \multicolumn{1}{c|}{Alg.~\ref{algorithm-basic}} &
          \multicolumn{1}{c|}{Alg.~\ref{algorithm-regular-*}} &
          \multicolumn{1}{c|}{Alg.~\ref{algorithm-motzkin}} \\
      \hline 
      0  & 1  &&& \\ 
      1  & 2  &&& \\ 
      2  & 7  &&& \\ 
      3  & 31 &&& \\ 
      4  & 153  &&& \\ 
      5  & 834  &&& \\ \hline
      6  & 4839  &&& \\ 
      7  & 29\ 612 &&&  \\ 
      8  & 188\ 695 &3&& \\  
      9  & 1243\ 746 &30 & 2 & \\ 
      10 & 8428\ 597 & - & 2& \\ \hline
      11 & 58\ 476\ 481 &- & 12& \\ 
      12 & 413\ 893\ 789 &- & 81& \\ 
      13 & 2980\ 489\ 256 &- &640 & 2\\ 
      14 & 21\ 787\ 216\ 989 & -& 5424 & 18\\ 
      15 & 161\ 374\ 041\ 945 &- & 46\ 330 & 212\\ \hline
      16 & 1209\ 258\ 743\ 839 &- & - & 1917\\ 
      17 & 9155\ 914\ 963\ 702 &- & - & 16\ 200\\ 
      18 & 69\ 969\ 663\ 242\ 487 &- & - & 136\ 980\\ 
      19 & 539\ 189\ 056\ 700\ 627 &- & - & 1096\ 320 \\
      \hline
    \end{tabular}
  \end{center}
  \vspace{-5mm}
  \caption{The number of idempotents in the Motzkin monoid $\M_n$, and the
  time in seconds to calculate these numbers using
  Algorithms~\ref{algorithm-basic},~\ref{algorithm-regular-*}
  and~\ref{algorithm-motzkin}.}
 \label{tab-EMn}
\end{table}

Tables \ref{tab-ErJn}, \ref{tab-ErKn} and \ref{tab-ErMn} give values of $|E_r(\J_n)|$, $|E_r(\K_n)|$ and $|E_r(\M_n)|$, respectively, for values of $n\leq10$; recall that $E_r(S)$ is the set of all idempotents of $S$ of rank $r$, where $S$ is any of $\J_n$, $\K_n$ or $\M_n$; cf.~Sequences A281441, A281442~and A269736 on \cite{OEIS}.  These values were calculated using the Semigroups package for GAP \cite{GAP}.  Higher values of these sequences could be calculated, by modifying Algorithms \ref{algorithm-jones}, \ref{algorithm-kauffman} and \ref{algorithm-motzkin} in light of Theorem \ref{thm-ErMJK}, but we have not done so.
Note also that for odd $n$, $|E_1(\K_n)|$ is a \emph{meandric number}; see Sequences A005315 and A005316 in \cite{OEIS}, and also \cite{DFGG1997,LaCroix2003}.

\begin{table}[ht]%
\begin{center}
\begin{tabular}{|c|rrrrrrrrrrr|}
\hline
$n\setminus r$	&0	&1	&2	&3	&4	&5	&6	&7	&8	&9	&10 \\
\hline
\phantom{1}0 & 1 &&&&&&&&&&\\
\phantom{1}1 & &1 &&&&&&&&&\\
\phantom{1}2 & 1&& 1 &&&&&&&&\\
\phantom{1}3 & &4&& 1 &&&&&&&\\
\phantom{1}4 & 4&& 7&& 1 &&&&&&\\
\phantom{1}5 & &25&& 10&& 1 &&&&&\\
\phantom{1}6 & 25&& 57&& 13&& 1 &&&&\\
\phantom{1}7 & &196&& 98&& 16&& 1 &&&\\
\phantom{1}8 & 196&& 522&& 148&& 19&& 1 &&\\
\phantom{1}9 & &1764&& 1006&& 207&& 22&& 1 &\\
10 & 1764&& 5206&& 1673&& 275&& 25&& 1 \\
\hline
\end{tabular}
\end{center}
\vspace{-5mm}
\caption{The number of rank $r$ idempotents in the Jones monoid $\J_n$.}
\label{tab-ErJn}
\end{table}


\begin{table} 
\begin{center}
\begin{tabular}{|c|rrrrrrrrrrr|}
\hline
$n\setminus r$	&0	&1	&2	&3	&4	&5	&6	&7	&8	&9	&10 \\
\hline
\phantom{1}0 & 1 &&&&&&&&&&\\
\phantom{1}1 & &1 &&&&&&&&&\\
\phantom{1}2 & 0&& 1 &&&&&&&&\\
\phantom{1}3 & &2&& 1&&&&&&&\\
\phantom{1}4 & 0&& 4&& 1 &&&&&&\\
\phantom{1}5 & &8&& 6&& 1 &&&&&\\
\phantom{1}6 & 0&& 22&& 8&& 1 &&&&\\
\phantom{1}7 & &42&& 40&& 10&& 1 &&&\\
\phantom{1}8 & 0&& 140&& 62&& 12&& 1 &&\\
\phantom{1}9 & &262&& 288&& 88&& 14&& 1 &\\
10 & 0&& 992&& 492&& 118&& 16&& 1 \\
\hline
\end{tabular}
\end{center}
\vspace{-5mm}
\caption{The number of rank $r$ idempotents in the Kauffman monoid $\K_n$.}
\label{tab-ErKn}
\end{table}


\begin{table}[ht]%
\begin{center}
\begin{tabular}{|c|rrrrrrrrrrr|}
\hline
$n\setminus r$	&0	&1	&2	&3	&4	&5	&6	&7	&8	&9	&10 \\
\hline
\phantom{1}0 & 1 &&&&&&&&&&\\
\phantom{1}1 & 1& 1 &&&&&&&&&\\
\phantom{1}2 & 4& 2& 1 &&&&&&&&\\
\phantom{1}3 & 16& 11& 3& 1 &&&&&&&\\
\phantom{1}4 & 81& 48& 19& 4& 1 &&&&&&\\
\phantom{1}5 & 441& 266& 93& 28& 5& 1 &&&&&\\
\phantom{1}6 & 2601& 1492& 549& 152& 38& 6& 1 &&&&\\
\phantom{1}7 & 16\ 129& 9042& 3211& 947& 226& 49& 7& 1 &&&\\
\phantom{1}8 & 104\ 329& 56\ 712& 20\ 004& 5784& 1480& 316& 61& 8& 1 &&\\
\phantom{1}9 & 697\ 225& 369\ 689& 127\ 676& 37\ 048& 9432& 2169& 423& 74& 9& 1 &\\
10 & 4787\ 344& 2477\ 806& 841\ 945& 241\ 268& 62\ 149& 14\ 402& 3036& 548& 88& 10& 1   \\
\hline
\end{tabular}
\end{center}
\vspace{-5mm}
\caption{The number of rank $r$ idempotents in the Motzkin monoid $\M_n$.}
\label{tab-ErMn}
\end{table}


As noted earlier, even though the monoid $\PP_n$ of all planar partitions of degree $n$ is isomorphic to the Jones monoid $\J_{2n}$ of degree $2n$, this is not true of their twisted versions, $\PP_n^\tau$ and $\J_{2n}^\tau=\K_{2n}$.  In general,~$\J_{2n}$ contains more idempotents than $\PP_n^\tau$.  The methods of this paper do not lead to algorithms for counting the idempotents of $\PP_n^\tau$.  However, for completeness, we used GAP \cite{GAP} to calculate the number of these idempotents for $n\leq10$.  Table \ref{tab-EPPn} gives the total number of these idempotents, while Table \ref{tab-ErPPn} gives the number of idempotents of a fixed rank; cf.~Sequences A286867~and A289620~on \cite{OEIS}.



%
%

\begin{table} 
  \begin{center}
    \begin{tabular}{|c|rrrrrrrrrrr|}
    \hline
$n$ &0&1&2&3&4&5&6&7&8&9&10\\
\hline
$|E(\PP_n^\tau)|$ & 1&1&6&44& 362& 3226& 30\ 488& 301\ 460 & 3090\ 020 & 32\ 618\ 046 & 345\ 515\ 557 \\
\hline
    \end{tabular}
  \end{center}
\vspace{-5mm}
  \caption{The number of idempotents in the twisted planar partition monoid $\PP_n^\tau$.}
  \label{tab-EPPn}
\end{table}
  


\begin{table}[b]%
\begin{center}
{\small
\begin{tabular}{|c|rrrrrrrrrrr|}
\hline
$n\setminus r$	&0	&1	&2	&3	&4	&5	&6	&7	&8	&9	&10 \\
\hline
\phantom{1}0& 1&&&&&&&&&&\\
\phantom{1}1& 0& 1&&&&&&&&&\\
\phantom{1}2& 0& 5& 1&&&&&&&&\\
\phantom{1}3& 0& 33& 10& 1 &&&&&&&\\
\phantom{1}4& 0& 253& 93& 15& 1 &&&&&&\\
\phantom{1}5& 0& 2147& 880& 178& 20& 1 &&&&&\\
\phantom{1}6& 0& 19\ 593& 8599& 1982& 288& 25& 1 &&&&\\
\phantom{1}7& 0& 188\ 837& 86\ 762& 21\ 723& 3684& 423& 30& 1 &&&\\
\phantom{1}8  &0&1899\ 107&900\ 997&238\ 419&44\ 767&6111&583&35&1&&\\
\phantom{1}9  &0&19\ 761\ 209&9595\ 264&2638\ 114&531\ 656&81\ 606&9388&768&40&1&\\
10 &0&211\ 447\ 863&104\ 447\ 385&29\ 503\ 900&6255\ 952&1044\ 248&136\ 740&13\ 640&978&45&1\\
\hline
\end{tabular}
}
\end{center}
\vspace{-5mm}
  \caption{The number of rank $r$ idempotents in the twisted planar partition monoid $\PP_n^\tau$.}
\label{tab-ErPPn}
\end{table}

\subsection*{Acknowledgements}

We would like to thank Dr Jan De Beule (Vrije Universiteit Brussel) and the
Department of Mathematics (WE01) of Ghent University for the use of their IBM
power8.

\footnotesize
\def\bibspacing{-1.1pt}
\bibliography{biblio}{}
\bibliographystyle{plain}

\end{document}